\documentclass[12pt,a4paper,fleqn]{article}
\pdfoutput=1

\usepackage[english]{babel}
\RequirePackage{ucs}
\RequirePackage[utf8x]{inputenc}
\RequirePackage[T1]{fontenc}
\RequirePackage{amsthm,amsmath}
\RequirePackage{amssymb,amstext}
\RequirePackage{amsfonts,stmaryrd}
\RequirePackage{amsbsy} 
\RequirePackage{mathrsfs}
\RequirePackage{bm}
\RequirePackage{bbm}
\RequirePackage{aliascnt,ifthen}
\RequirePackage{graphicx}
\RequirePackage{color,calc}
\RequirePackage{times}
\RequirePackage{booktabs}
\RequirePackage{dsfont} 
\usepackage{url}
\usepackage{verbatim}
\usepackage{subcaption}
\usepackage[pdfpagemode=UseOutlines ,plainpages=false,
hyperindex=true,colorlinks=true]{hyperref}
	\makeatletter
	\Hy@breaklinkstrue
	\makeatother

\definecolor{darkred}{rgb}{0.6,0,0.1}
\definecolor{darkgreen}{rgb}{0,0.5,0}
\definecolor{darkblue}{rgb}{0,0,0.5}

\hypersetup{colorlinks ,linkcolor=gray ,filecolor=darkgreen, urlcolor=darkblue ,citecolor=black}

\RequirePackage{paralist}
\RequirePackage{relsize}
\usepackage{enumitem}

\usepackage{natbib}
\usepackage{hypernat}
\renewcommand{\cite}{\citet}
\usepackage{SCORDEMME-abrev-package}

\definecolor{dgreen}{rgb}{0,0.5,0}
\definecolor{dblue}{rgb}{0,0,0.5}
\definecolor{dred}{rgb}{0.6,0.0,0.1}
\definecolor{dgold}{rgb}{0.5,0.3,0.0}
\definecolor{dvio}{rgb}{0.6,0.3,0.5}
\definecolor{gray}{rgb}{0.5,0.5,0.5}
\definecolor{dbraun}{rgb}{.5,0.2,0}

\newcommand{\dgrau}{\color{gray}}

\newcommand{\colre}{dred}
\newcommand{\colas}{dblue}
\newcommand{\colrem}{dgold}
\newcommand{\colil}{dgreen}

\newtheoremstyle{styre}% name
  {1.1\topsep}%      Space above
  {\topsep}%      Space below
  {\normalfont\itshape}%         Body font
  {}%         Indent amount (empty = no indent, \parindent = para indent)
  {\color{\colre}}% Thm head font
  {.}%        Punctuation after thm head
  {.5em}%     Space after thm head: " " = normal interword space;
  {\thmname{\textbf{#1}\xspace}{\xspace\dgrau\thmnumber{#2}}\thmnote{\xspace\textit{\small(#3)}}}%         Thm head spec (can be left empty, meaning `normal')

\newtheoremstyle{styas}% name
  {1.1\topsep}%      Space above
  {\topsep}%      Space below
  {\normalfont\itshape}%         Body font
  {}%         Indent amount (empty = no indent, \parindent = para indent)
  {\color{\colas}}% Thm head font
  {.}%        Punctuation after thm head
  {.5em}%     Space after thm head: " " = normal interword space;
  {}%         Thm head spec (can be left empty, meaning `normal')

\newtheoremstyle{styrem}% name
  {1.1\topsep}%      Space above
  {\topsep}%      Space below
  {\normalfont\itshape}%         Body font
  {}%         Indent amount (empty = no indent, \parindent = para indent)
  {\color{\colrem}}% Thm head font
  {.}%        Punctuation after thm head
  {.5em}%     Space after thm head: " " = normal interword space;
  {\thmname{\textbf{#1}\xspace}{\xspace\dgrau\thmnumber{#2}}\thmnote{\xspace\textit{\small(#3)}}}%{}%         Thm head spec (can be left empty, meaning `normal')

\newtheoremstyle{styil}% name
  {1.1\topsep}%      Space above
  {\topsep}%      Space below
  {\normalfont\rmfamily}%         Body font
  {}%         Indent amount (empty = no indent, \parindent = para indent)
  {\color{\colil}}% Thm head font
  {.}%        Punctuation after thm head
  {.5em}%     Space after thm head: " " = normal interword space;
  {\thmname{\textbf{#1}\xspace}{\xspace\dgrau\thmnumber{#2}}\thmnote{\xspace\textit{\small(#3)}}}%         Thm head spec (can be left empty, meaning `normal')

\newtheoremstyle{stypro}%
	{0.5\topsep}%space above
	{1.1\topsep}%space below
	{\upshape}%		body font
	{}%				indent amount
	{}%	theorem head font
	{.}%			punctuation after theorem head
	{.5em}%			space after theorem head
	{\thmnote{\textit{#3}}}%         Thm head spec (can be left empty, meaning `normal')

\theoremstyle{styre}\newtheorem{pr}{Proposition}[section]
\newaliascnt{co}{pr}
\theoremstyle{styre}\newtheorem{co}[co]{Corollary}
\aliascntresetthe{co}
\newaliascnt{thm}{pr}
\theoremstyle{styre}\newtheorem{thm}[thm]{Theorem}
\aliascntresetthe{thm}
\newaliascnt{lem}{pr}
\theoremstyle{styre}\newtheorem{lem}[lem]{Lemma}
\aliascntresetthe{lem}
\newaliascnt{rem}{pr}
\theoremstyle{styrem}\newtheorem{rem}[rem]{Remark}
\aliascntresetthe{rem}
\newaliascnt{il}{pr}
\theoremstyle{styil}\newtheorem{il}[il]{Illustration}
\aliascntresetthe{il}
\theoremstyle{styas}

\theoremstyle{styas}

\theoremstyle{stypro}\newtheorem*{pro}{}

\newcommand{\remEnd}{{\scriptsize\textcolor{\colrem}{\qed}}}
\newcommand{\proEnd}{{\scriptsize\textcolor{\colre}{\qed}}}

\usepackage{cleveref}
\crefname{pr}{\color{\colre}Proposition}{\color{\colre}Propositions}
\crefname{co}{\color{\colre}Corollary}{\color{\colre}Corollaries}
\crefname{thm}{\color{\colre}Theorem}{\color{\colre}Theorems}
\crefname{lem}{\color{\colre}Lemma}{\color{\colre}Lemmata}
\crefname{ass}{\color{\colas}Assumption}{\color{\colas}Assumptions}
\crefname{de}{\color{\colas}Definition}{\color{\colas}Definitions}
\crefname{rem}{\color{\colrem}Remark}{\color{\colrm}Remarks}
\crefname{il}{\color{\colil}Illustration}{\color{\colil}Illustrations}

\usepackage{environ}
\NewEnviron{te}{%
{\BODY}%
}

\numberwithin{equation}{section} 

\makeatletter
\newcommand{\mylabel}[2]{#2\def\@currentlabel{#2}\label{#1}}
\makeatother

\newcommand{\setListe}[5][3ex]{\setlength{\itemsep}{#2}\setlength{\topsep}{#3}\setlength{\leftmargin}{#4}\setlength{\rightmargin}{#5}\setlength{\labelwidth}{#1}}
\newcommand{\setListeStandard}{\setListe{0ex}{.5ex}{4ex}{0ex}}
\newcounter{ListeN}
\renewcommand{\theListeN}{(\roman{ListeN})}

\newenvironment{resListeN}[1][~]%
{\setcounter{ListeN}{0}\renewcommand{\theListeN}{\normalfont\rmfamily\color{\colre}(\roman{ListeN})}\begin{list}{\theListeN}%
{\usecounter{ListeN}\setListeStandard #1}}
{\end{list}}

{\begin{list}{}%
{\setListeStandard #1}}
{\end{list}}

\makeatletter% 
\def\@fnsymbol#1{\ensuremath{\ifcase#1\or * \or ** \or 2 \or 3 \or  *\or  \star \or 4\or  , \or 
g\or h\or i\else\@ctrerr\fi}}% 
\makeatother% 

\author{\begin{minipage}{.45\textwidth}\center{\sc Sergio Brenner Miguel}\;\thanks{Institut f\"ur Angewandte
    Mathematik, M$\Lambda$THEM$\Lambda$TIKON, Im Neuenheimer Feld 205,
  D-69120 Heidelberg, Germany, e-mail:
  \url{{brennermiguel|johannes}@math.uni-heidelberg.de}}\\\small Ruprecht-Karls-Universit\"{a}t Heidelberg\\\null
\end{minipage} \and \begin{minipage}{.45\textwidth}\center{\sc 
  Fabienne Comte} \thanks{CNRS, MAP5 UMR 8145,
  F-75006 Paris, France, e-mail:
  \url{fabienne.comte@parisdescartes.fr}}\\\small Universit\'e de Paris\\\null\end{minipage}\and\begin{minipage}{.45\textwidth}\center{\sc 
  Jan Johannes}$\;^*$\\\small Ruprecht-Karls-Universit\"{a}t Heidelberg\\\null\end{minipage}}
\date{\today} 
\title{Spectral cut-off regularisation for density estimation under multiplicative measurement errors} 
\begin{document} 
\maketitle 
% --------------------------------------------------------------------
% <<Abstract>>
% --------------------------------------------------------------------
\begin{abstract}
  We study the non-parametric estimation of an unknown density $f$ with support on $\pRz$ based on an i.i.d. sample with multiplicative measurement errors. 
  The proposed fully-data driven procedure is based on the estimation of the Mellin
  transform of the density $f$, a regularisation of the inverse of
  the Mellin transform by a spectral cut-off  and a
  data-driven model selection  in order to deal with the
  upcoming bias-variance trade-off. We introduce and discuss further
  \textit{Mellin-Sobolev spaces} which characterize the
  regularity of the unknown density $f$ through the decay of its Mellin
  transform. Additionally,  we show minimax-optimality over \textit{Mellin-Sobolev spaces} of the
 data-driven density estimator and hence its adaptivity. 
\end{abstract} 
{\footnotesize
\begin{tabbing} 
\>minimax theory, inverse problem, adaptation\\[.2ex] 
\end{tabbing}}

% --------------------------------------------------------------------
% <<Content>>
% --------------------------------------------------------------------
%======================================================================================================================
%                                                                 
% Title: Introduction
% Author: Sergio Brenner Miguel and Jan JOHANNES, IAM, Ruprecht-Karls Universität Heidelberg, Deutschland  
% 
% Date: %%ts latex start%%[2020-07-30 Thu 19:10]%%ts latex end%%
%
% ======================================================================================================================
% ====================================================================
% Section <<Intro>>
% ====================================================================
\section{Introduction}\label{i}
% ....................................................................
% Motivation, Objectiv
% ....................................................................
\begin{te} 
  In this paper we are interested in estimating the unknown density
  $f:\pRz \rightarrow \pRz$ of
  a positive random variable $X$ given  independent and identically
  distributed (i.i.d.) copies of $Y=XU$, where $X$ and $U$ are
  independent of each other and $U$ has a known density $g:\pRz
  \rightarrow \pRz$. In this setting the density  $f_Y:\pRz
  \rightarrow \pRz$ of $Y$ is given by
  \begin{equation*}
    f_Y(y)=[f * g](y):= \int_{0}^{\infty} f(x)g(y/x) x^{-1}dx\quad\forall y\in\pRz
  \end{equation*}
  such that $*$ denotes multiplicative convolution. The estimation of
  $f$ using an i.i.d. sample $Y_1, \dots, Y_n$ from $f_Y$ is thus an
  inverse problem called
  multiplicative deconvolution.  \\
  \cite{Vardi1989} and \cite{VardiZhang1992} introduce and study
  intensively \textit{multiplicative censoring}, which corresponds to
  the particular multiplicative deconvolution problem with
  multiplicative error $U$ uniformly distributed on $[0,1]$.  This
  model is often applied in survival analysis as explained and
  motivated in \cite{Van-EsKlaassenOudshoorn2000}. The estimation of
  the cumulative distribution function of $X$ is considered in
  \cite{VardiZhang1992} and
  \cite{AsgharianCaroneFakoorothers2012}. Series expansion methods are
  studied in \cite{AndersenHansen2001} treating the model as an
  inverse problem. The density estimation in a multiplicative
  censoring model is considered in \cite{BrunelComteGenon-Catalot2016}
  using a kernel estimator and a convolution power kernel
  estimator. Assuming a uniform error distribution on an interval
  $[1-\alpha, 1+\alpha]$ for $\alpha\in (0,1)$ \cite{ComteDion2016}
  analyze a projection density estimator with respect to the Laguerre
  basis.  \cite{BelomestnyComteGenon-Catalot2016} study a
  beta-distributed error $U$.\\
  In this work, covering all those three variations of the
  multiplicative censoring model we consider a density estimator using
  the Mellin transform and a spectral cut-off regularization of its
  inverse, which borrows ideas from
  \cite{BelomestnyGoldenshlugerothers2020}.  The key to the analysis
  of the multiplicative deconvolution problem is the multiplication
  theorem of the Mellin transform $\Melop$, which roughly states
  $\Melop[f_Y]=\Melop[f]\Melop[g]$ for a density $f_Y=f
  *g$. Exploiting the multiplication theorem
  \cite{BelomestnyGoldenshlugerothers2020} introduce a kernel density
  estimator of $f$ allowing more generally $X$ and $U$ to be
  real-valued. Moreover, they point out that the following widely used
  naive approach is a special case of their estimation
  strategy. Transforming the data by applying the logarithm the model
  $Y=XU$ writes $\log(Y)=\log(X)+\log(U)$. In other words,
  multiplicative convolution becomes convolution for the
  $\log$-transformed data. As a consequence, the density of $\log(X)$
  is eventually estimated employing usual strategies for
  non-parametric deconvolution problems (see for example
  \cite{Meister2009}) and then transformed back to an estimator of $f$.
  However, it is difficult to interpret regularity conditions on the
  density of $\log(X)$. Furthermore, the analysis of a global risk of
  an estimator using this naive approach is challenging as
  \cite{ComteDion2016}
  pointed out.\\
  Our strategy differs in the following way.  Making use of the
  multiplication theorem of the Mellin transform and applying an
  additional spectral cut-off on the inversion of the Mellin-transform we
  define a density estimator. We measure the accuracy of the estimator
  by introducing a global risk in terms of a weighted
  $\LpA$-norm. Exploiting properties of the Mellin transform we
  characterize the underlying inverse problem and natural regularity
  conditions which borrow ideas from the inverse problems community
  (\cite{EnglHanke-BourgeoisNeubauer2000}).  The regularity conditions expressed
  in the form of \textit{Mellin-Sobolev spaces} and their relations to
  the analytical properties of the density $f$ are discussed in more
  details. The proposed estimator, however, involves a tuning
  parameter which is selected by a data-driven method. We establish an
  oracle inequality for the fully-data driven spectral cut-off
  estimator under fairly mild assumptions on the error density
  $g$. 
  Moreover we show that uniformly over \textit{Mellin-Sobolev spaces} the proposed
  data-driven estimator is minimax-optimal. Precisely, we state both an upper
  bound for the mean weighted integrated squared error of the fully-data driven spectral
  cut-off  estimator and a general lower bound for estimating the
  density $f$  based on i.i.d. copies from $f_Y=f *g$.
\end{te}

% ....................................................................
% Organisation
% ....................................................................

The paper is organized as follows. In \cref{ag} we collect properties
of the Mellin transform. We explain our general estimation strategy by
first introducing and analyzing an estimator based on direct
observations $X_1,\dots, X_n$ from $f$. The estimator relies on an
inversion of the Mellin-transform which we stabilize using a spectral
cut-off.  Exploiting then the multiplication theorem of the
Mellin-transform we propose a fully-data driven estimator of $f$ based
on the sample $Y_1,\dots,Y_n$. We derive an oracle type upper bound
for its mean weighted integrated squared error. We finish the
theoretical part by showing in \cref{mm} that our fully-data driven
estimator is minimax optimal over \textit{Mellin-Sobolev spaces} for a
large class of error densities $g$. Finally, results of a simulation
study are reported in \cref{si} which visualize the reasonable finite
sample performance of our estimators. The proof of \cref{ag} and
\cref{mm} are postponed to the appendix.

%%% Local Variables:
%%% mode: latex
%%% TeX-master: "_0SCORDEMME"
%%% End:

%======================================================================================================================
%                                                                 
% Title: Data-driven procedure
% Author: Sergio Brenner Miguel and Jan JOHANNES, IAM, Ruprecht-Karls Universität Heidelberg, Deutschland  
% 
% Date: %%ts latex start%%[2020-07-30 Thu 19:49]%%ts latex end%%
%
% ======================================================================================================================

\section{Adaptive weighted $\LpA$ estimation}\label{ag}
In this section we introduce the Mellin transform and collect some of its properties.  
For a more detailed introduction we refer to \cite{ParisKaminski2001} and \cite{BarucqMattesiTordeux2015}.
\paragraph{Mellin transform}
Let $\Lz^{1,\text{loc}}_{\pRz}$ denote the set of locally integrable real-valued functions. For  $h\in \Lz^{1,\text{loc}}_{\pRz}$ the Mellin transform of $h$ in the point $c+it\in \Cz$ is defined by 
\begin{align}\label{eq:mel;def}
\Mela h c(t) := \mathcal M[h](c+it):=\int_0^{\infty} x^{c-1+it} h(x)dx
\end{align}
provided that the integral is absolutely convergent.  If there exists a $c\in \Rz$ such that the mapping $x\mapsto x^{c-1}h(x)$ is integrable over $\pRz$ then the region $\Xi_h \subseteq \Cz$ of absolute convergence of the integral in \eqref{eq:mel;def} is either a vertical strip $\{s+it\in \Cz: s\in(a,b), t\in \Rz \}$ for $a<b$ with $c\in(a,b)$ or a vertical line $\{c+it\in \Cz: t\in \Rz\}$. In the case that $\Xi_h$ is a vertical strip the function $s+it\mapsto\Mela h s(t)$ is analytical on this strip. In the literature $\Xi_h$ is often called \textit{strip of analyticity}. 
In the following illustration we give some techniques to determine $\Xi_h$.
\begin{il}\label{il:mel:well}
	Note that for any density $h_1\in \Lz_{\pRz}^{1,\text{loc}}$  the vertical strip $\{1+it: t\in \Rz\}$ belongs to $\Xi_{h_1}$, and hence the Mellin transform $\Mel[h_1]$ is well-defined. Furthermore 
	the region $\Xi_{h_2}$ of $h_2\in \Lz_{\pRz}^{1,\text{loc}}$ is a superset of the vertical strip $\{c+it: c\in (a,b), t\in \Rz\}$ if $h_2(x) = \mathrm O(x^{-a+\varepsilon})$ for $x \downarrow 0$ and $h_2(x)= \mathrm O(x^{-b-\varepsilon})$ for $x\rightarrow \infty$ and for all $\varepsilon>0$ small enough.  
	The latter property implies that for the family $(g_k)_{k\in \Nz}$ with $g_k(x):= \1_{(0,1)}(x) k(1-x)^{k-1}$ for a $k\in \Nz$ and $x\in \pRz$ the function $\Melop[g_k]$ is analytical on $\{c+it: c>0, t\in \Rz\}$ since for all $ b, x>0, \varepsilon \in (0,b)$   holds $|x^{\varepsilon}g_k(x)|\vee|x^{b-\varepsilon} g_k(x)| \leq k$. 
\end{il}
For $ c \in \Xi_h $  the inversion formula of the Mellin transform
is given by 
\begin{align}\label{eq:Mel:inv}
h(x)= \frac{1}{2\pi } \int_{-\infty}^{\infty} x^{-c-it} \Mela h c(t) dt,\text{ for }x>0,
\end{align} 
whenever $\Xi_h$ is not a vertical line
(c.f. \cite{ParisKaminski2001}) or alternatively if the function $t\mapsto  \Mela h c(t) dt$ is square integrable over $\Rz$ (c.f. \cite{BarucqMattesiTordeux2015}). 
Considering \cref{il:mel:well} we see that the assumption on $\Xi_h$
not being a vertical line is rather weak. It is fulfilled for almost
all functions in the upcoming theory. However, in all the other cases
we make use  of the second assumption without further reference.\\
It can be shown that $\int_0^{\infty} h^2(x)x^{2\alpha-1} dx< \infty$ for $\alpha\geq0$ implies that $\alpha \in \Xi_h$ and also that
$\int_0^{\infty} h^2(x)x^{2\alpha-1} dx = (2\pi)^{-1}
\int_{-\infty}^{\infty} |\Mela h\alpha(t)|^2dt$. This result combined
with the Mellin inversion formula implies an isometry in the following
way. For $\alpha\geq 0$ define the weight function
$\omega_{\alpha}(x):=x^{2\alpha-1}, x\in \Rz,$ and the corresponding
weighted norm   by $\|h\|_{\omega_{\alpha}}^2 := \int_0^{\infty}
h^2(x) \omega_{\alpha}(x)dx $ for a measurable function. Denote by
$\LpA(\omega_{\alpha})$ the set of all measurable functions with
finite $\|\, .\,\|_{\omega_{\alpha}}$-norm and by $\langle h_1, h_2
\rangle_{\omega_{\alpha}} := \int_0^{\infty}  h_1(x)
h_2(x)\omega_{\alpha}(x)dx$ for $h_1, h_2\in \LpA(\omega_{\alpha})$
the corresponding weighted scalar product. Similarly, define $\Lp[2]\Rz:=\{ h:\Rz \rightarrow \Cz\, \text{ measurable }: \|h\|_{\Rz}^2:= (2\pi)^{-1} \int_{-\infty}^{\infty} h(t)\overline{h(t)} dt <\infty \}.$ Now,  following \cite{BarucqMattesiTordeux2015} both operators $\Melop_{\alpha}: \LpA(\omega_{\alpha}) \rightarrow \Lp[2] \Rz$ and 
\begin{align}\label{eq:Mel:inv:op}
\Melop_{\alpha}^{-1}: \Lp[2] \Rz \rightarrow
  \LpA(\omega_{\alpha}),\quad  h \mapsto (x\mapsto \Melop_{\alpha}^{-1}[h](x):=(2\pi)^{-1} \int_{-\infty}^{\infty} x^{-\alpha-it} h(t) dt)
\end{align} are isometries. We will first construct an estimator for $f$ given an i.i.d. sample $X_1,\dots, X_n$, that is the direct case, and afterwards we construct  an estimator based on the i.i.d. sample $Y_1,\dots,Y_n$, which constitutes the censored case.

\paragraph{Case of direct observation} 
In this paragraph we define the estimator of $f\in \mathbb L^2_{\pRz}(\omega_{\alpha})$ given the direct observations $X_1, \dots, X_n$ by using the Mellin transform and spectral cut-off regularised inverse. 
Since $f\in \LpA(\omega_{\alpha})$ the Mellin transform $\Mela {f} \alpha$ is well-defined and a natural unbiased estimator of $\Mela {f} \alpha[t]$ for each $t\in \Rz$ is given by $\widehat{\mathcal M}_{ \alpha}(t):= n^{-1} \sum_{j=1}^n X_j^{\alpha-1+it}$. It is worth pointing out that this estimator is bounded in $t\in \Rz$ by  $|\widehat{\mathcal M_{\alpha}}(t)| \leq |\widehat{\mathcal M}_{\alpha}(0)|$ which is finite almost surely. Thus $\1_{[-k,k]} \widehat{\mathcal M_{\alpha}}\in \Lz_{\Rz}^2$ for all $k\in \pRz$ which allows us to apply the operator in \eqref{eq:Mel:inv:op} to define 
\begin{align}\label{eq:est1}
\widehat f_k(x):= \Melop_{\alpha}^{-1}[\1_{[-k,k]}\widehat{\mathcal M}_{\alpha}](x) = \frac{1}{2\pi} \int_{-k}^k x^{-\alpha-it} \widehat{\mathcal M}_{ \alpha}(t) dt,\, \text{for } \, x\in \pRz,
\end{align}
as an unbiased estimator of $f_k:=\Melop_{\alpha}^{-1}[\1_{[-k,k]} \Melop_{\alpha}[f]]$.
Additionally, we see that $\|f-f_k\|_{\omega_{\alpha}}^2 = \pi^{-1} \int_{k}^{\infty} |\Mela f{\alpha}(t)|^2 dt$ tending to zero for $k$ going to infinity, that is $f_k$ approximates $f$ in the weighted $\LpA$ sense. \\
Furthermore, we have $\|\widehat f_k \|_{\omega_{\alpha}}^2 = (2\pi)^{-1} \int_{-k}^k |\widehat{\mathcal M}_{ \alpha}(t)|^2 dt$ which as a random variable has only a finite first moment if and only if the by $f$ induced distribution has a finite $2(\alpha-1)$ moment by application of Fubini-Tonelli theorem. We state the following proposition which implies the consistency of the estimator for a suitable choice of the cut-off parameter $k\in \pRz$.  

\begin{pr}\label{dd:pr:con}
	Assume that $f\in \LpA(\omega_{\alpha})$ and that $\sigma^2:=\E_f(X^{2(\alpha-1)}) < \infty$. Then, for all $k\in \Rz$ we have
	\begin{align*}
	\E_{f}^n(\|f- \widehat f_k\|_{\omega_{\alpha}}^2) 
	& \leq\|f-f_k \|_{\omega_{\alpha}}^2  +  \pi^{-1}\sigma^2 k n^{-1}.
	\end{align*}
	By choosing $k=k_n$ such that $n^{-1} k_n \rightarrow 0$ and $k_n \rightarrow \infty$, $\widehat f_{k_n}$ is a consistent estimator of $f$.
\end{pr}
The proof of \cref{dd:pr:con} can be found in \cref{a:ag}.
Note that, if one would like to consider the case $\alpha=1/2$, that is $\omega_{\alpha}=1$, it is necessary to assume a finite first  moment of $X^{-1}$. On the other hand, a less restrictive situation would be to consider the case of $\alpha=1$ which needs no additional moment condition on $f$ (respectively on $g$) since in this case $\sigma^2=1$. The corresponding weight function would be $\omega(x):= \omega_1(x)=x$ for $x \in \pRz$. 
For the estimation of densities without a compact support assumption, the intervals far away from zero are of special interest. A weighted $\LpA$-norm could model this and may allow us to capture more interesting characteristics of the density like being heavy-tailed or compactly supported.
 From now on, we will restrict ourselves to this special case while we want to remark that the upcoming theory can be expanded to different values of $\alpha\geq0$ under additional assumptions. \\
 Before we start to define the estimator in the case of censored observation let us briefly discuss the upcoming bias term $\|f-f_k\|_{\omega_{\alpha}}^2$ in \cref{dd:pr:con}. A natural condition which allows a more sophisticated study on the bias, is to assume that the Mellin-transform decays with a polynomial rate, that is for $s\geq 0$,
 \begin{align}\label{in:eq:MelSob}
 f \in\Wz^s:=\{h \in \LpA: |h|^2_s:=\int_{-\infty}^{\infty} |\Mel[h](t)|^2 (1+t^2)^s dt< \infty\}.
 \end{align} 
 The definition of these spaces strongly resemble to the frequently considered Sobolev spaces which are defined as $W^s:=\{H \in \Lz_{\Rz}^2 : \int_{-\infty}^{\infty} |\mathcal F[H](t)|^2 (1+t^2)^{s} dt < \infty\}$ for $s\geq 0$ where $\mathcal F[H](t):=\int_{-\infty}^{\infty} H(x) \exp(-ixt)dx$ denotes the usual Fourier-transform for an integrable function $H:\Rz\rightarrow \Rz$. To distinguish between them, we refer $\Wz^s$ as \textit{Mellin-Sobolev space} and $W^s$ as \textit{Fourier-Sobolev space}. But not just the general motivation seems to be similar as we can easily see by the following relationship between Fourier transformation and Mellin transformation. For $h\in \LpA(\omega)$ we have $\Mela h 1 = \Mela {\omega h}{0} = \mathcal F[(\omega h) \circ \varphi]$ with $\varphi: \Rz \rightarrow \pRz, x\mapsto \exp(-x)$. Setting $H:= (\omega h) \circ \varphi$ we see that $h \in \Wz^s$ is equivalent to $H\in W^s$. Therefore, it does not seem to be a suprise that it is possible to characterise the \textit{Mellin-Sobolev spaces} via analytical properties as done in the case of the \textit{Fourier-Sobolev spaces.} This is stated in the following proposition while its proof is postponed to the appendix.
 \begin{pr}\label{mm:pr:reg}
 	Let $s \in \Nz$. Then $f \in \Wz^s$ if and only if $f$ is $s-1$-times continuously differentiable where $f^{(s-1)}$ is locally absolutely continuous with derivative $f^{(s)}$ and  $\omega^j f^{(j)} \in \LpA(\omega)$ for all $j\in \nset{0,s}$.
 \end{pr}
	If a positive random variable $R$ has the density $h$,  the real-valued random variable $T=\log(R)$ has the density $f_T(y)=H(-y)= (\omega h)\circ \varphi(-y)$ for $y \in \Rz$. Again, we observe the strong connection between the Mellin transform of $h$ and the Fourier transform of $H$.  This has the following interesting implication for the multiplicative measurement error model.

 \begin{rem}\label{mm:rem:dec}
	As already mentioned the application of the logarithm to the random variable $Y$,  $Z:=\log(Y)=\log(X)+\log(U)=: V + \varepsilon$, where $V$ and $\varepsilon$ are independent, can be used to transfer the model to a regular deconvolution setting. This technique was used for instance 
	 by \cite{ComteDion2016}. Given an estimator $\widehat f_V$ of
         $f_V$ it is possible to derive an estimator of $f$ through $\widehat f:=\omega^{-1} \widehat f_V\circ \log$. In fact, one can show that $\| \widehat f_V- f_V\|_{\Rz}^2 = \| \widehat f - f\|_{\omega}^2$ which illustrate the difficulties which arise when considering the global risk. Furthermore, in the deconvolution approach via a Fourier transformation one usually assume that the densities $f_V, f_{\varepsilon}\in \Lz_{\Rz}^2$ which again correspond to the fact that $f, g\in \LpA(\omega)$ using the considerations above. 
\end{rem}
In the next part we define an estimator of $f$ based on the censored observation $Y_1, \dots, Y_n$. Since $Y=XU$, In the multiplicative measurement error model we would  need to assume that both $X$ and $U$ have a finite $-1$-moment to consider the unweighted norm. Especially the latter scenario would exclude several interesting examples, for instance the multiplicative censoring model where $U$ is uniformly distributed on $[0,1]$.

\paragraph{Case of censored observation } 
The key property which makes the Mellin transform useful for multiplicative deconvolution is that for two function $h_1,h_2\in \Lz_{\pRz}^{1,\text{loc}}$ with $c \in \Xi_{h_1}\cap \Xi_{h_2}$, 
\begin{align}\label{eq:mel:multi}
\Mela {h_1*h_2}{c} (t) = \Mela{h_1}{c}(t)\cdot\Mela{h_2} c (t) \quad \text{ for } t\in \Rz. 
\end{align} We will refer to it from now on as the multiplication theorem. In the context of deconvolution a similar property adressing the convolution and its Fourier transform is frequently used to construct deconvolution estimators. Since $f$ and $g$ are both densities we have $1 \in \Xi_f \cap \Xi_g$ which implies that for all $t\in \Rz$, $\Mel[f_Y](t)= \Mel[f](t) \Mel[g](t)$.
Under the assumption  that  $\Mel[g](t)\neq 0$ for all $t\in \Rz$, which we do in the upcoming theory without further reference, we see that using \eqref{eq:mel:multi} we can express the functions $(f_k)_{k\in\pRz}$ in the following way
\begin{align*}
	 f_k(x)= \frac{1}{2\pi} \int_{-k}^{k} x^{-1-it} \Mel[f](-t)dt =\frac{1}{2\pi} \int_{-k}^{k} x^{-1-it} \frac{\Mel[f_Y](t)}{\Mel[g](t)}dt
\end{align*} 
for $k,x \in \pRz$. Similar to the direct case we define our estimator by replacing $\Mel[f_Y](t)$ with its empirical counterpart $\widehat{\mathcal M}(t) := \widehat{\mathcal M}_1(t)=\frac{1}{n} \sum_{j=1}^n Y_j^{it}$ to define the estimator 
\begin{align}\label{eq:est2}
\widehat f_k(x):= \frac{1}{2\pi} \int_{-k}^k  x^{-1-it} \frac{\widehat{\M}(t)}{\Mel[g](t)} dt \quad \text{ for } \, x,k>0.
\end{align}
The following mild assumption on the error density ensures that the estimator is well-defined, for all $k\in \pRz$,
\begin{align*}
\int_{-k}^k |\Mel[g](t)|^{-2} dt< \infty %\text{ and }|\Mel[g](t)|>0  \text{ for all } |t|\leq k
  \tag{\textbf{[G0]}}.
\end{align*}
Note that  $\Mel[\widehat f_k](t)= \1_{[-k,k]} (t)\frac{\widehat
  \M(t)}{\Mel[g](t)}$ by definition of $\widehat f_k$, and hence the
estimator defined in \eqref{eq:est2} and  \eqref{eq:est1} coincide
when setting $\Mel[g](t)=1$ for $t\in \Rz$. The proof  of the next
proposition is  very similar to the proof of \cref{dd:pr:con} and thus omitted.
\begin{pr}\label{dd:pr:con2}
	Assume that $f\in \LpA(\omega)$ and that \textbf{[G0]} holds. Then for all $k\in \Nz$ we have
	\begin{align*}
	\E_{f_Y}^n(\|f- \widehat f_k\|_{\omega}^2) 
	& \leq \|f-f_k \|_{\omega}^2  + (2\pi n)^{-1} \int_{-k}^k |\Mel[g](t)|^{-2} dt.
	\end{align*}
	By choosing $k=k_n$ such that $n^{-1} \int_{-k_n}^{k_n} |\Mel[g](t)|^{-2} dt \rightarrow 0$ and $k_n \rightarrow \infty$, $\widehat f_{k_n}$ is a consistent estimator of $f$.
\end{pr}
Let us now have a closer look at the second summand in \cref{dd:pr:con2} which bounds the variance term of the estimator. In the following, we use for two functions $f,g$  the notation
$f\sim g$ over a set $A\subset \Rz$ if the function $f/g$ is
bounded away both from zero and infinity over the set $A$. In analogy to the usual deconvolution setting and to the work of \cite{BelomestnyGoldenshlugerothers2020} we say that the error density is \textit{smooth} if there exist parameters $\gamma, \tau_1\in \pRz$ such that
\begin{align*}
\tag{\textbf{[G1]}}  \forall |t| \geq \tau_1:  | \Mel[g](\tau )| \sim t^{-\gamma} \text{ and }
\forall |t| \leq \tau_1:  |\Mel[g](t)| \sim 1.
\end{align*} 
Now \textbf{[G1]} implies that $\int_{-k}^k |\Mel[g](t)|^{-2} dt \leq C_g k^{2\gamma+1}$ where $C_g>0$ is a positive constant introduced in the following corollary.

\begin{co}\label{dd:co:con2}
	Assume that $f\in \LpA(\omega)$ and that \textbf{[G1]} holds. Then for all $k\in \Nz$ we have
	\begin{align*}
	\E_{f_Y}^n(\|f- \widehat f_k\|_{\omega}^2) 
	& \leq \|f-f_k \|_{\omega}^2  + C_g(2\pi n)^{-1} k^{2\gamma+1}
	\end{align*}
	where $C_g$ is a constant only dependent on $g$. By choosing $k_n$ such that $n^{-1} k_n^{2\gamma+1}\rightarrow 0$ and $k_n \rightarrow \infty$, $\widehat f_{k_n}$ is a consistent estimator of $f$.
\end{co}
Under the assumptions of \cref{dd:co:con2} it is natural to restrict the set of suitable parameters $k$ to $\mathcal K_n:=\nset{1, K_n}$ with $K_n:=n^{1/(2\gamma+1)}$ and to choose $k_n:\in\argmin\{\|f-f_k\|_{\omega}^2 + C_g(2\pi n)^{-1} k^{2\gamma+1} : k \in \mathcal K_n\}$. Unfortunately, this choice is not feasible since it depends on the unknown density $f$ itself. We note that the bias $\|f-f_k\|_{\omega}^2 = \|f\|_{\omega}^2 - \|f_k \|_{\omega}^2$ behaves like  $-\|f_k\|_{\omega}^2$. Exchanging $-\|f_k\|_{\omega}^2$ with its empirical counterpart $-\|\widehat f_k\|_{\omega}^2$ we define a fully data-driven model selection $\widehat k$  by   
\begin{align}\label{eq:data:driven}
\widehat k \in \argmin \{-\|\widehat f_k\|_{\omega}^2 + \mathrm{pen}(k) : k \in \mathcal K_n\} \quad \text{where} \quad\mathrm{pen}(k):= \chi k^{2\gamma+1} n^{-1}
\end{align}
for $\chi>0.$
The following theorem shows that this procedure is adaptive up to a negligeable term.
\begin{thm}\label{dd:thm:ada}
Assume that $f\in  \LpA(\omega)$, \textbf{[G1]} and that $\|\omega f_Y\|_{\infty} := \sup_{y>0} |y f_Y(y)| <\infty$. Then for $\chi > 12 C_g \pi^{-1}$ 
\begin{align*}
\E_{f_Y}^n (	\| f- \widehat f_{\widehat k} \|_{\omega}^2) \leq 4\inf_{k\in\mathcal K_n}\big(\|f-f_k\|_{\omega}^2 +\mathrm{pen}(k) \big) + C(\|\omega f_Y\|_{\infty},g) n^{-1}
\end{align*}
where $C(\|\omega f_Y\|_{\infty},g)>0$ is a constant depending on $\|\omega f_Y\|_{\infty}$ and $g$.
\end{thm}

The proof of  \cref{dd:thm:ada} is postponed to appendix \cref{a:ag}.
The assumption that $\|\omega f_Y\|_{\infty}<\infty$ is rather weak. In fact, since $1\in \Xi_f \cap \Xi_g$ we are able to write $|yf_Y(y)|= |y\frac{1}{2\pi} \int_{-\infty}^{\infty} y^{-1-it}\Mel[f](t)\Mel[g](t) dt| \leq \|f\|_{\omega} (\int_{-\infty}^{\infty} |\Mel[g](t)|^2 dt)^{1/2} < \infty$ if $\gamma >1/2$ in \textbf{[G1]}. \\
The last assertion establishes an oracle inequality assuming a smooth
error density as in \textbf{[G1]}. For a \textit{super smooth} error
density with exponentially decay of its  Mellin transform (see
\cite{BelomestnyGoldenshlugerothers2020}) a result similar to
\cref{dd:thm:ada} can be derived from \cref{dd:lem:talapply} in the
appendix provided the upper bound $K_n$ and the penalty terms
$\mathrm{pen}(k)$, $k\in\nset{1, K_n}$  are
choosen accordingly. However, we omit the details, since the minimax theory presented in the next
chapter does not cover a \textit{super smooth} error density.

%%% Local Variables:
%%% mode: latex
%%% TeX-master: "_0SCORDEMME"
%%% End:

%======================================================================================================================
%                                                                 
% Title: Minimax theory
% Author: Sergio Brenner Miguel and Jan JOHANNES, IAM, Ruprecht-Karls Universität Heidelberg, Deutschland  
% 
% Date: %%ts latex start%%[2020-01-28 Tue 15:41]%%ts latex end%%
%
% ======================================================================================================================
\section{Minimax theory}\label{mm}

In this section we develop the minimax theory for the proposed estimator in \cref{ag}. Over the \textit{Mellin-Sobolev} spaces we derive an upper and lower bound for the mean weighted integrated squared error, which are equal up to a multiplicative constant, showing that our estimator is minimax-optimal over these spaces.

\paragraph{Regularity assumptions}
	Let us define for  $s\geq 0$ and the ellipsoids $\Wz^s(L):=\{h\in \Wz^s: |h|^2_s \leq L\}$ for any $L\geq 0$ which correspond to the \textit{Mellin-Sobolov spaces} defined in \eqref{in:eq:MelSob}.  We see that for any $f\in \Wz^s(L)$ we have $\int_{[-k,k]^c} |\Mel[f](t)|^2 dt \leq L k^{-2s}$ and $\|f\|_{\omega}^2= (2\pi)^{-1} \int_{-\infty}^{\infty} |\Mel[f](t)|^2 dt \leq L(2\pi)^{-1}$.  
 We denote the subset of densities  by 
 \begin{equation}\label{equation:dens_sobol}
 \rwcSobD{\wSob,\rSob}:= \{f\in \Wz^s(L): f \text{ is
 	a density}\}.
 \end{equation}
 Again, assuming \textbf{[G1]} implies that $\int_{-k}^k |\Mel[g](t)|^{-2} dt \leq C_g k^{2\gamma+1}$ where $C_g>0$ is a constant only dependent on the error density $g$. These considerations imply the following theorem whose proof is omitted.
 
 \begin{thm}\label{mm:thm:upper}
 	Assume that $g$ satisfies assumption \textbf{[G1]}. Then for $k_o:=n^{1/(2s+2\gamma+1)}$, 
 	\begin{align*}  
 		\sup_{f\in\rwcSobD{\wSob,\rSob}}\E_{f_Y}^n(\|f- \widehat f_{k_o}\|_{\omega}^2) \leq C(g, L, s) n^{-2s/(2s+2\gamma+1)}.
 	\end{align*}
 \end{thm}
As mentioned before for $\gamma>1/2$ we have $\|\omega f_Y\|_{\infty} \leq \|f\|_{\omega} \|g\|_{\omega}$. Thus, we can state the following corollary which is a direct consequence of \cref{dd:thm:ada} and \cref{mm:thm:upper}.
\begin{co}
	Assume that \textbf{[G1]} holds for $\gamma>1/2$. Then for $\chi >12 C_g\pi^{-1}$, 
	\begin{align*}
	\sup_{f\in \rwcSobD{\wSob, \rSob}}\E_{f_Y}^n (	\| f- \widehat f_{\widehat k} \|_{\omega}^2) \leq C(g,L,s) n^{-2s/(2s+2\gamma+1)},
	\end{align*}
	where $\chi$ is defined in \eqref{eq:data:driven}.
\end{co}

\begin{rem}
	For $g(x)= k(1-x)^{k-1} \1_{(0,1)}(x)$ with $k\in \Nz$,$x>0$ and $c>0$ its Mellin transform is given by $\Mela g c(t) =\prod_{j=1}^{k} \frac{j}{c-1+j+it} $ and satisfies \textbf{[G1]}  with $\gamma = k$. In fact this covers the model considered by  \cite{BelomestnyComteGenon-Catalot2016} as a generalisation of the multiplicative censoring where we consider a uniform distributed error density , that is $k=1$. Again, we can include the  case of direct observations, getting a rate of $n^{-2s/(2s+1)}$ over the ellipsoid.
\end{rem}
To prove that the rate of \cref{mm:thm:upper} is minimax-optimal over the ellipsoids under certain assumptions on $g$ we finish this section by stating a lower bound result. We want to emphasize that up to now we had no constraints on the support of $g$. To prove the lower bound we will need to assume that $g$ has a bounded support. For the sake of simplicity we will assume that $g$ has a support in $[0,1]$. We assume that there exist parameters $\gamma, \tau_1\in \pRz$ such that
\begin{align*}
\hspace*{-1.03cm}\tag{\textbf{[G1']}}  \forall |\tau| \geq \tau_1:  | \Mela g {1/2}(t)| \sim t^{-\gamma} \text{,}
\forall |\tau| \leq \tau_1:  |\Mela g {1/2}(t )| \sim 1 \text{ and } \forall x>1: g(x)=0.
\end{align*} 

\begin{thm}\label{theorem:lower_bound}
	Let $\wSob, \gamma \in\Nz$, assume that \textbf{[G1']} holds. Then there exist
	constants $\cst{g},L_{\wSob,g},n_{s,\gamma}>0$ such that for all
	$L\geq L_{\wSob,g}$, $n\geq n_{s,\gamma}$ and for any estimator $\hSo$ of $\So$ based
	on an i.i.d. sample $\Nsample{Y_j}$, 
	\begin{align*}
	\sup_{f\in\rwcSobD{\wSob,\rSob}}\E_{f_Y}^n(\|f- \widehat f\|_{\omega}^2) \geq \cst{g} n^{-2s/(2s+2\gamma+1)}.
	\end{align*}
\end{thm}
We want to emphasize that the error densities $(g_k)_{k\in \Nz}$ with $g_k(x)=k(1-x)^{k-1}\1_{(0,1)},x\in \pRz,$ fulfill both the assumption \textbf{[G1]} and \textbf{[G1']}. Thus, in this situation our estimation strategy is minimax-optimal. The proof of the lower bound can be extended to the case of directly observed $X_1,\dots,X_n$ or for different weight functions $\omega_{\alpha}, \alpha\geq 0$.

%%% Local Variables:
%%% mode: latex
%%% TeX-master: "_0SCORDEMME"
%%% End:

%======================================================================================================================
%                                                                 
% Title: Data-driven aggregation
% Author: Sergio Brenner Miguel and Jan JOHANNES, IAM, Ruprecht-Karls Universität Heidelberg, Deutschland  
% 
% Date: %%ts latex start%%[2020-01-22 Wed 15:02]%%ts latex end%%
%
% ======================================================================================================================
\section{Numerical study}\label{si}

Let us illustrate the performance of the estimator $\widehat f_{\widehat k}$ defined in \eqref{eq:est2} and \eqref{eq:data:driven} in the cases  $U\sim \mathrm{U}_{[0,1]}$, $U\sim \mathrm{U}_{(0.5,1.5)}$ and $U\sim \mathrm{Beta}_{(1,2)}$. For the density $g$ of an uniform distribution on $[0.5, 1.5]$ we get that $\Melop_{1}[g](t)= (1+it)^{-1} (1.5^{1+it}-0.5^{1+it}), t\in \Rz$, which corresponds to the case of $\gamma=1$ in \textbf{[G1]}.
We consider the densities
\begin{resListeN}
	\item\label{si:lag:i}  \textit{Gamma Distribution}: $f(x)=\frac{x^{4}}{4!}\exp(-x)$,
	\item\label{si:lag:ii} \textit{Gamma Mixture:} $f(x)=0.4 \cdot 3.2^2 x\exp(-3.2x) +0.6 \cdot \frac{6.8^{16} x^{15}}{15!}\exp(-6.8x)$,
	\item\label{si:lag:iii}  \textit{Beta Distribution}: $f(x)=\frac{1}{560} (0.5x)^3(1-0.5x)^4 \1_{[0,1]}(0.5x)$ and 
	\item\label{si:lag:iv} \textit{Weibull Distribution}: $f(x)= 2 x \exp(-x^2)$.
\end{resListeN}
By minimising an integrated weighted squared error over a family of histogram
densities with randomly  drawn partitions  and weights we select
$\chi= 1.2$, $\chi = 0.8$ and $\chi=0.01$ for the cases $\gamma=0,$ $\gamma=1$ and $\gamma = 2$, respectively, where $\chi$ is the penalty constant, see \eqref{eq:data:driven}.\\
In the direct case we compare the estimator $\widehat f_{\widehat k}$ with the data-driven density estimator $\widetilde f$ from the work of \cite{Brenner-MiguelJohannes2020} which is based on the adaptive aggregation of projection estimators with respect to the Laguerre basis. \\
\begin{minipage}{\textwidth}
\centering{\begin{minipage}[t]{0.32\textwidth}
		\includegraphics[width=\textwidth,height=35mm]{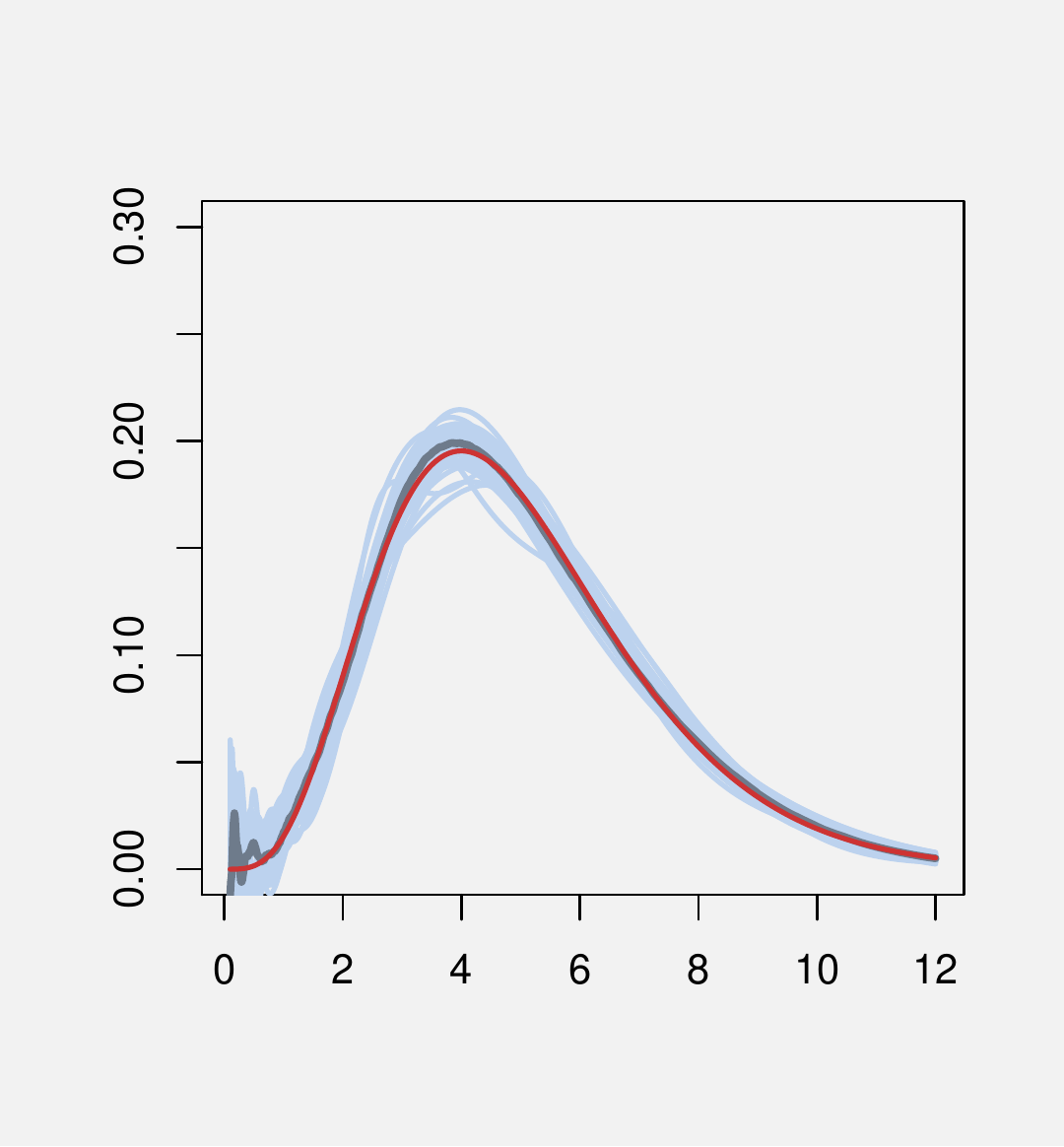}
	\end{minipage}
	\begin{minipage}[t]{0.32\textwidth}
		\includegraphics[width=\textwidth,height=35mm]{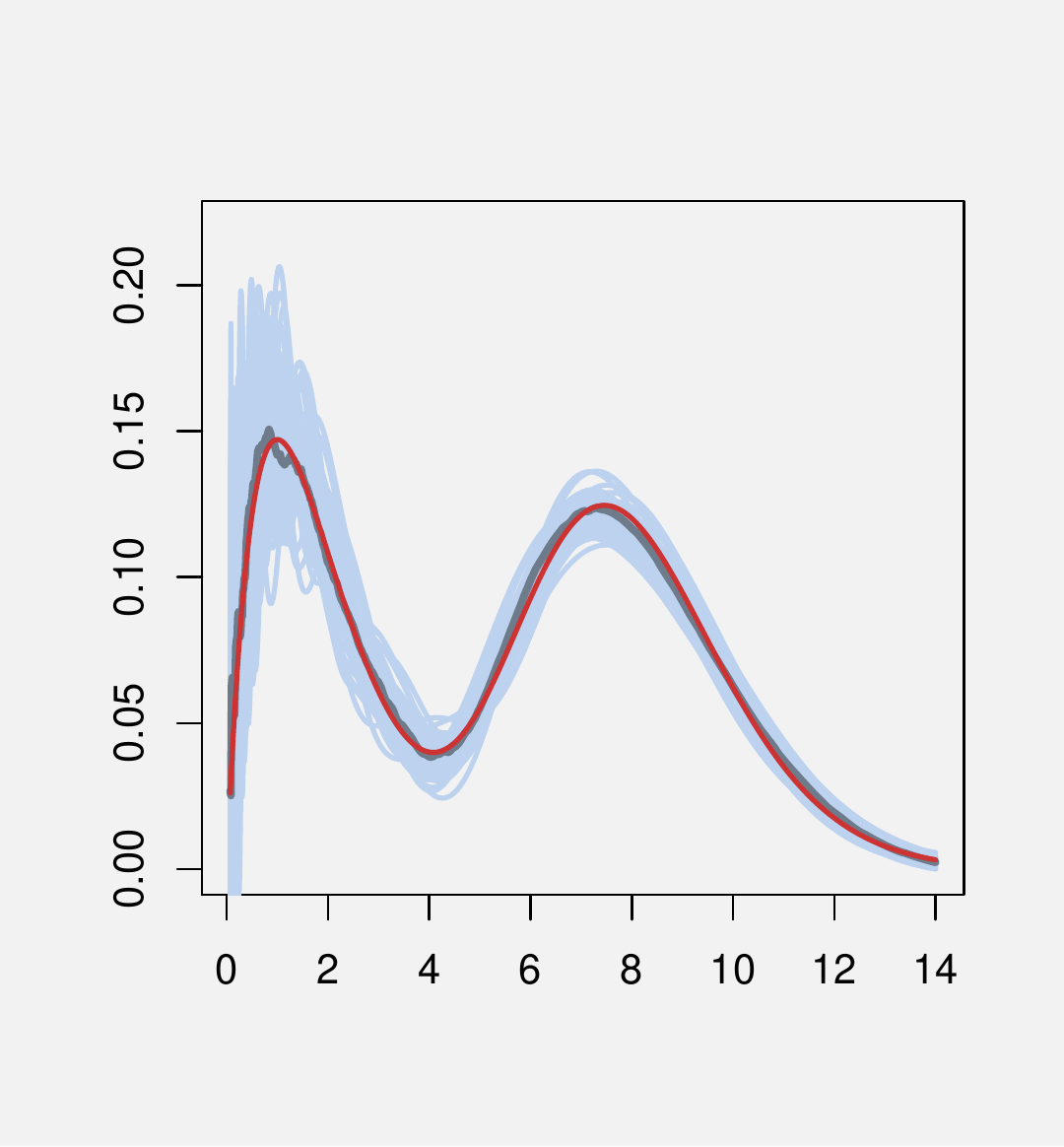}
	\end{minipage}
\begin{minipage}[t]{0.32\textwidth}
		\includegraphics[width=\textwidth,height=35mm]{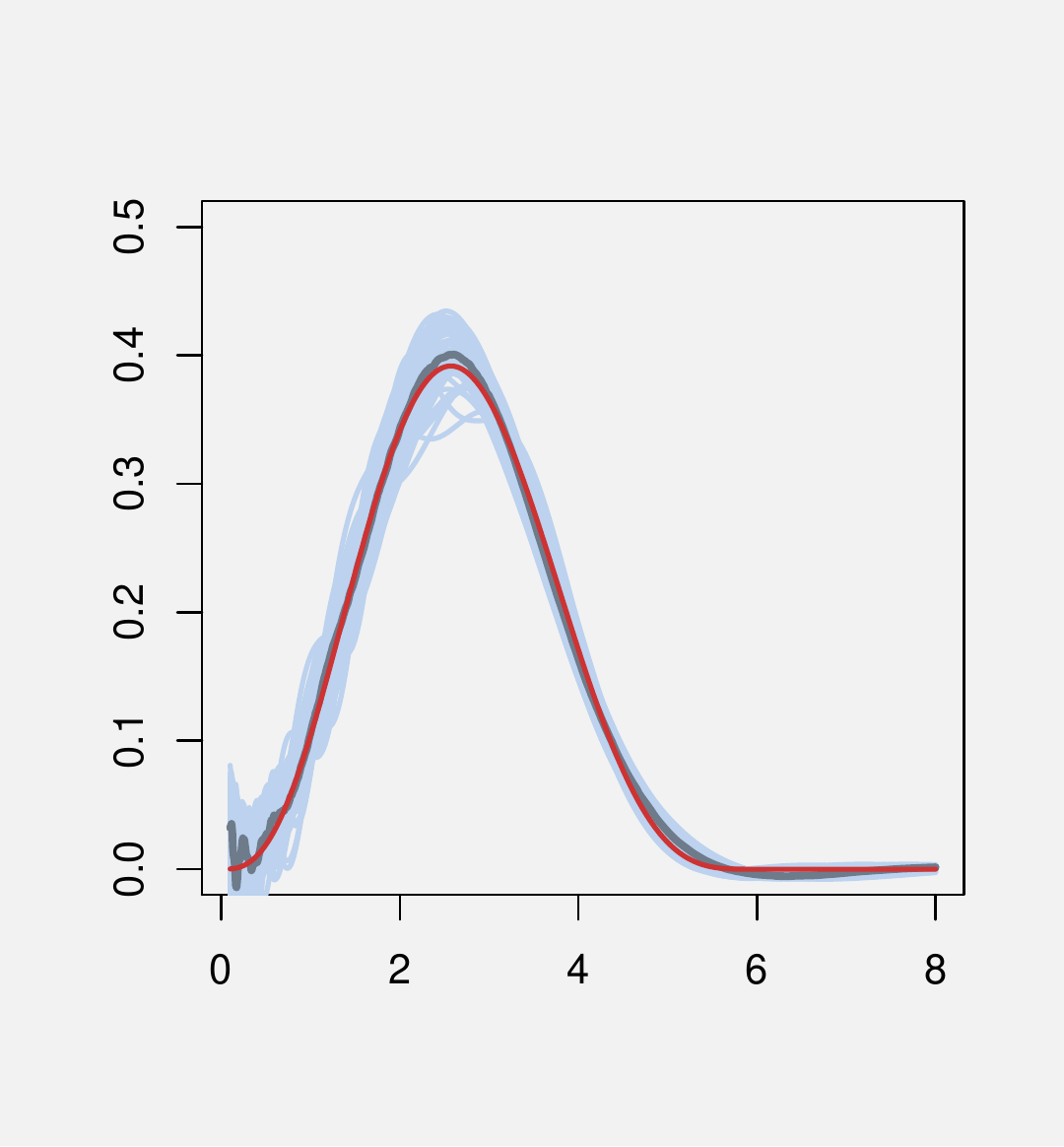}
	\end{minipage}}
\centering{\begin{minipage}[t]{0.32\textwidth}
		\includegraphics[width=\textwidth,height=35mm]{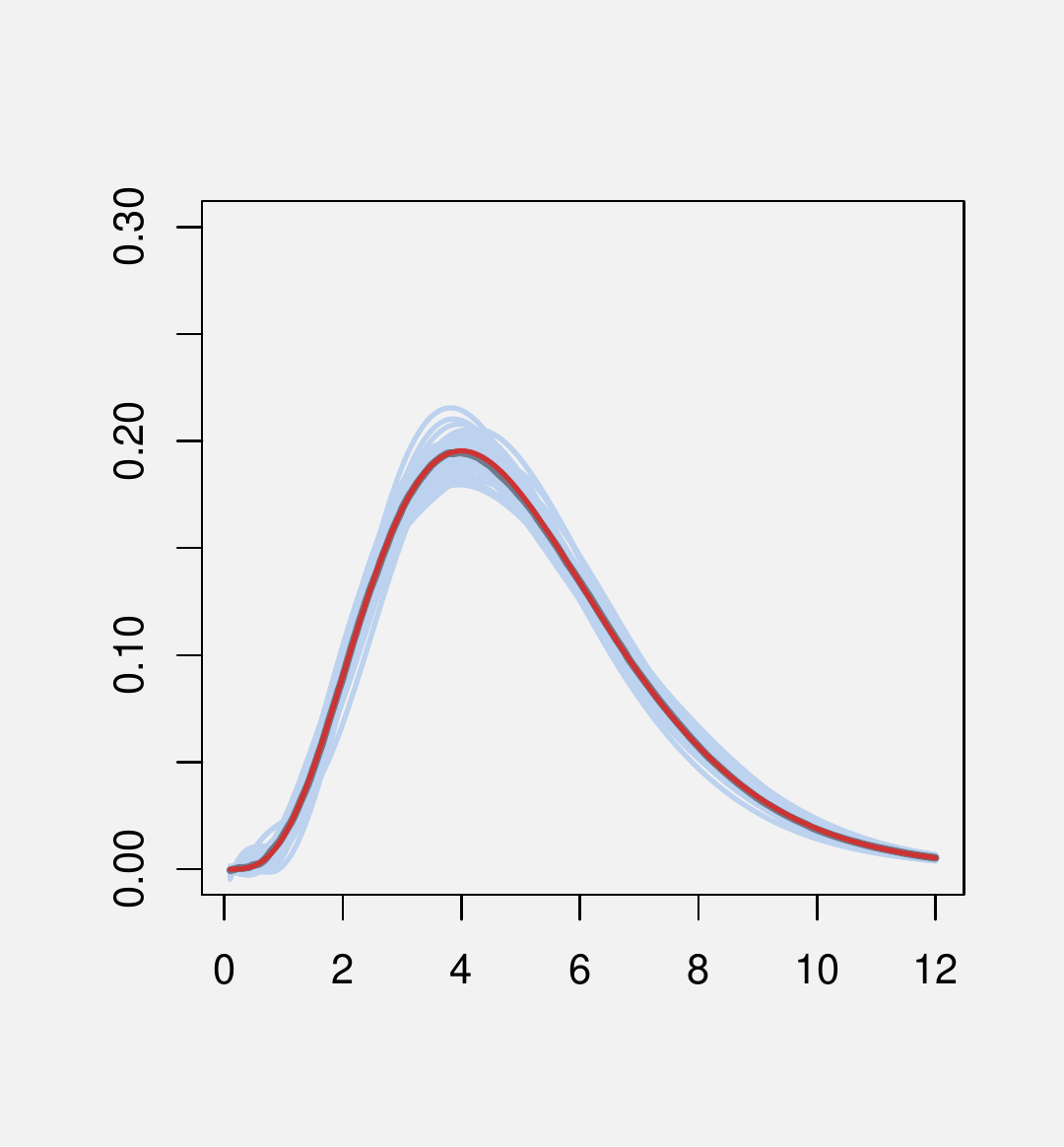}
	\end{minipage}
\begin{minipage}[t]{0.32\textwidth}
	\includegraphics[width=\textwidth,height=35mm]{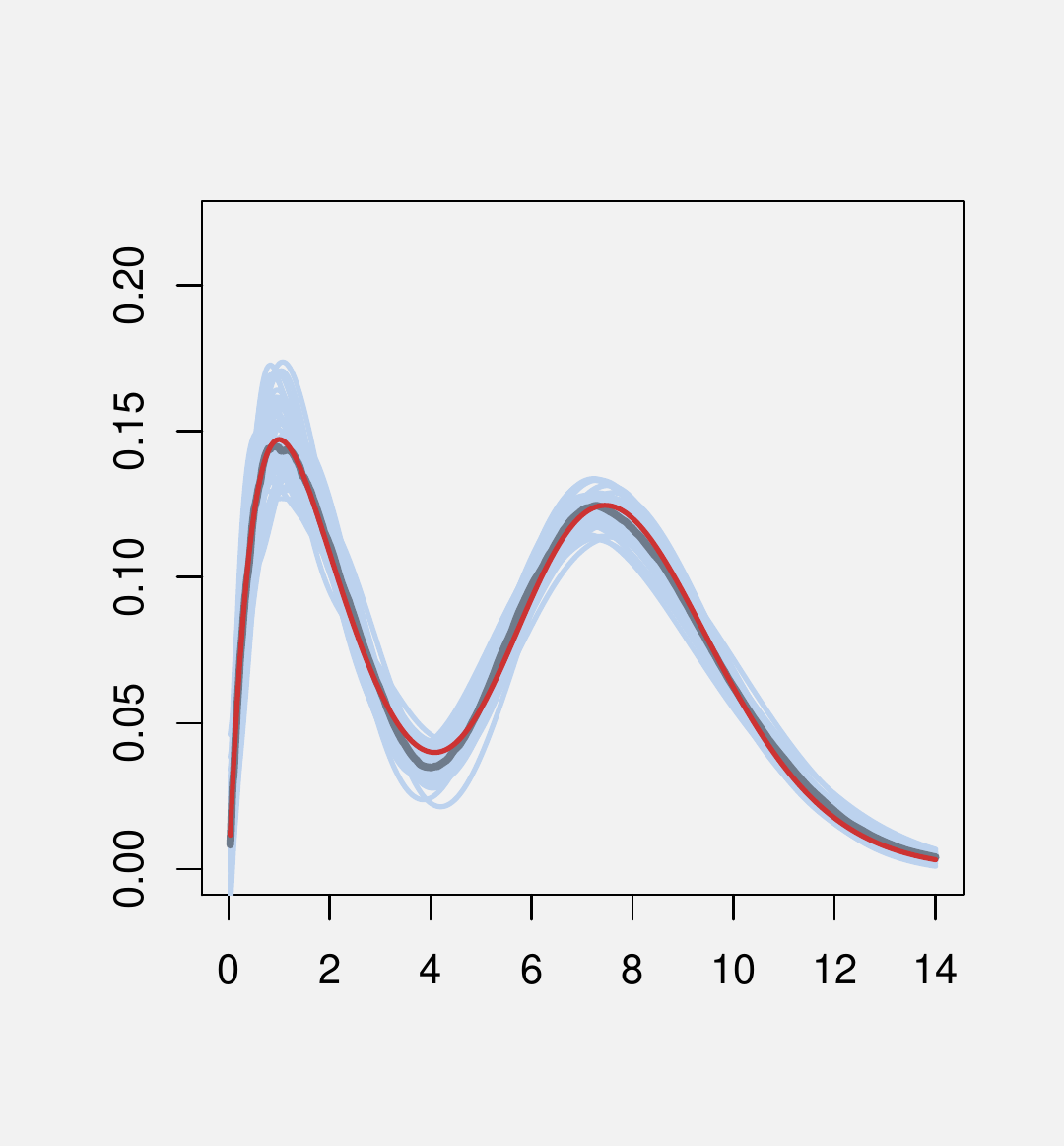}
\end{minipage}
\begin{minipage}[t]{0.32\textwidth}
	\includegraphics[width=\textwidth,height=35mm]{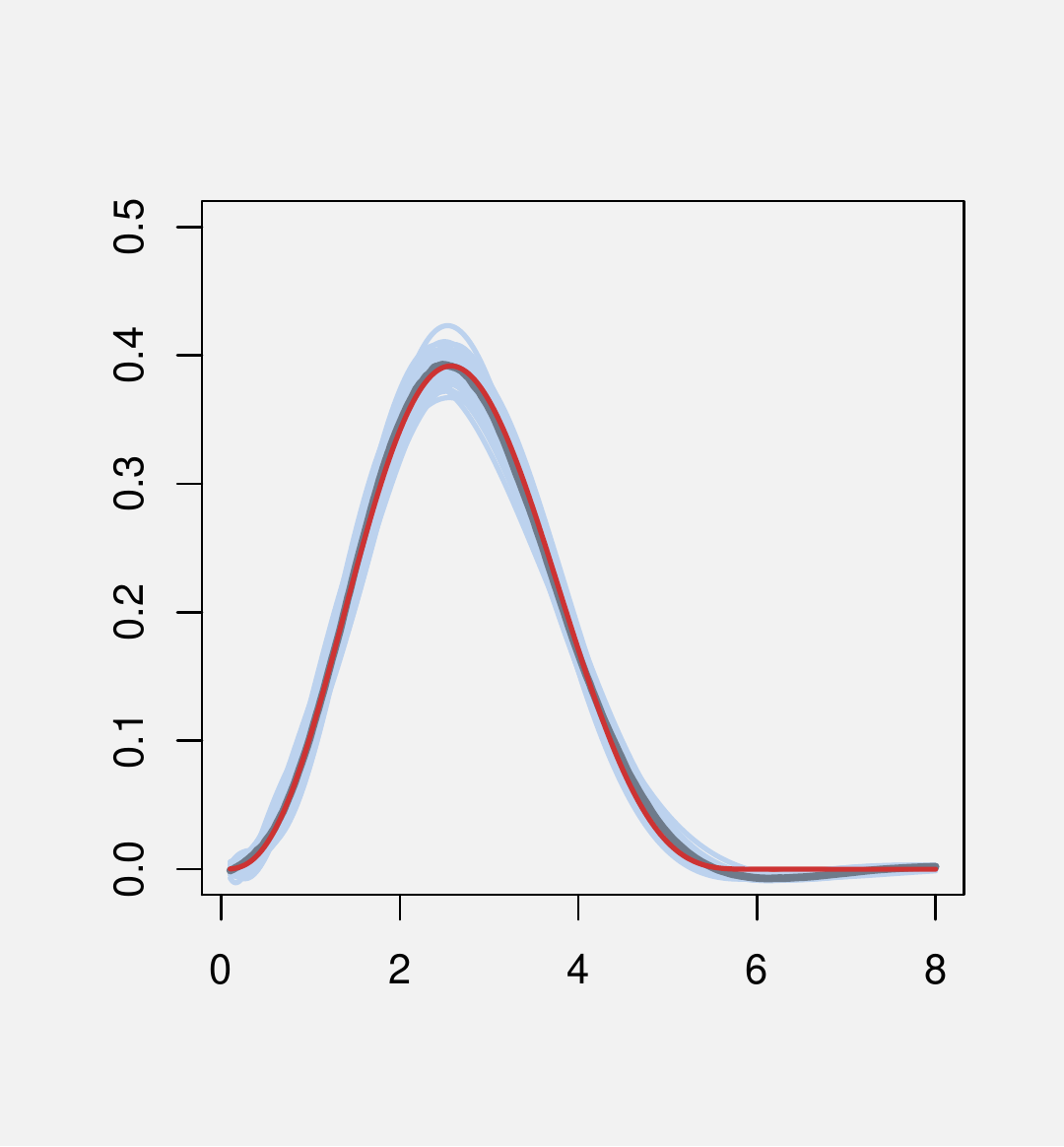}
\end{minipage}}
	
      \captionof{figure}{\label{figure:1}Considering 
        the estimators $\widehat f_{\widehat k}$ (top) and $\widetilde f$ (bottom) are depict for 
        50  Monte-Carlo simulations with  sample size $n=1000$ in the case \ref{si:lag:i} (left), \ref{si:lag:ii} (middle) and \ref{si:lag:iii} (right) with direct observations. The true density $\So$ is given by the red curve while the dark blue curve is the point-wise empirical median of the 50 estimates.}
\end{minipage}\\[2ex]
\begin{minipage}[t]{\textwidth}
\centerline{\begin{minipage}[t]{0.32\textwidth}
		\centering \includegraphics[width=\textwidth,height=35mm]{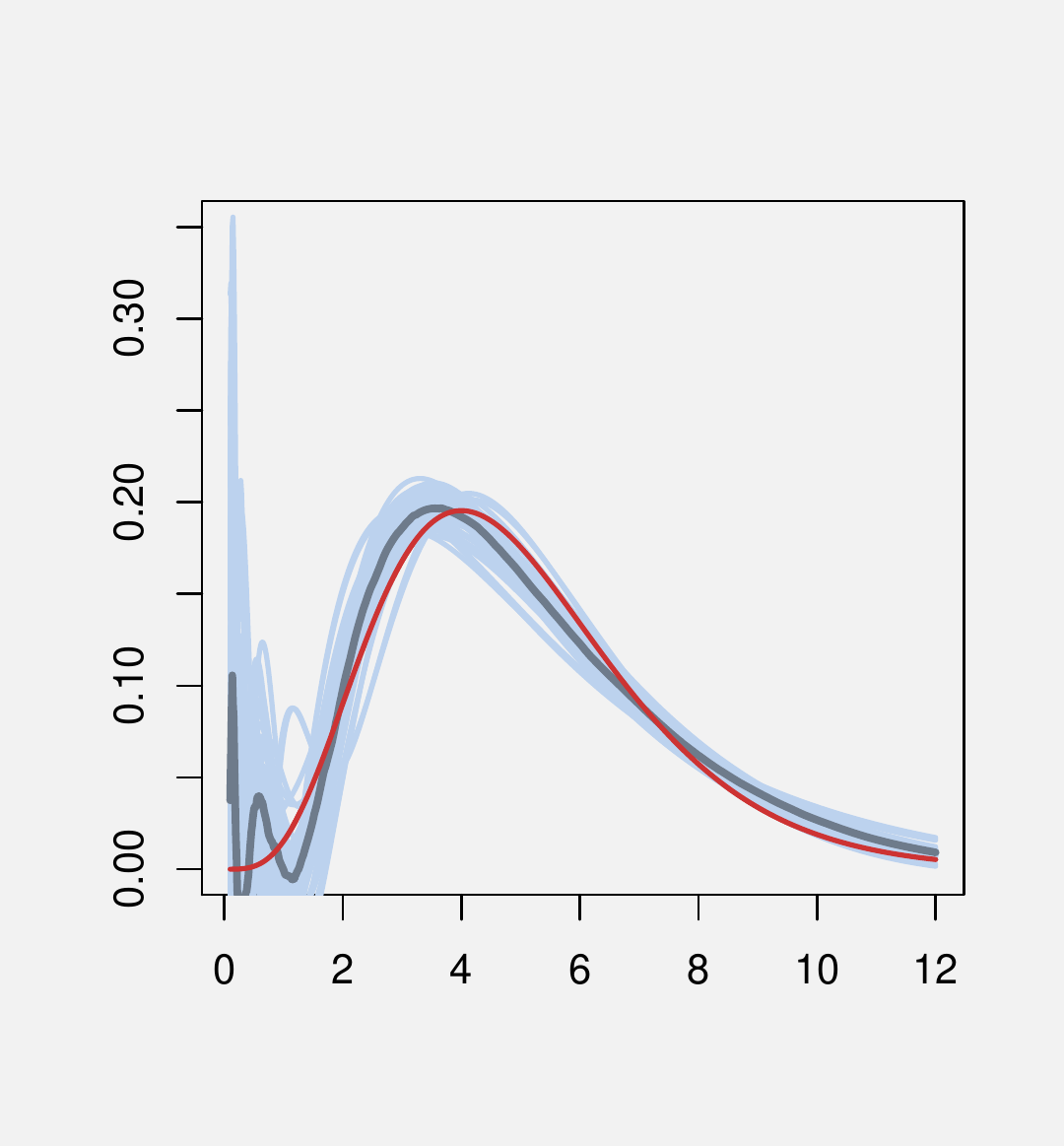}
	\end{minipage}
	\begin{minipage}[t]{0.32\textwidth}
	\centering	\includegraphics[width=\textwidth,height=35mm]{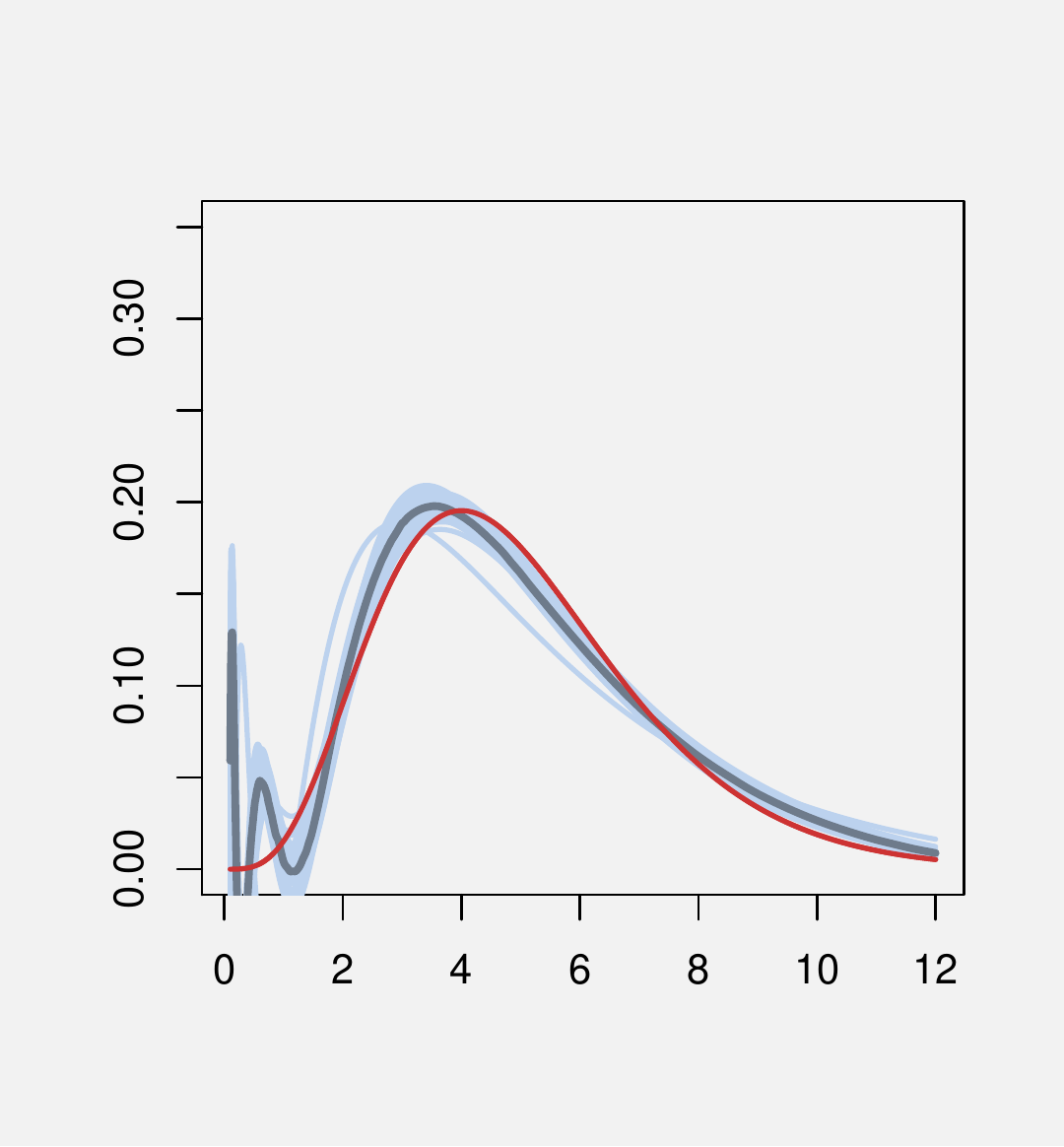}
	\end{minipage}
	\begin{minipage}[t]{0.32\textwidth}
		\centering	\includegraphics[width=\textwidth,height=35mm]{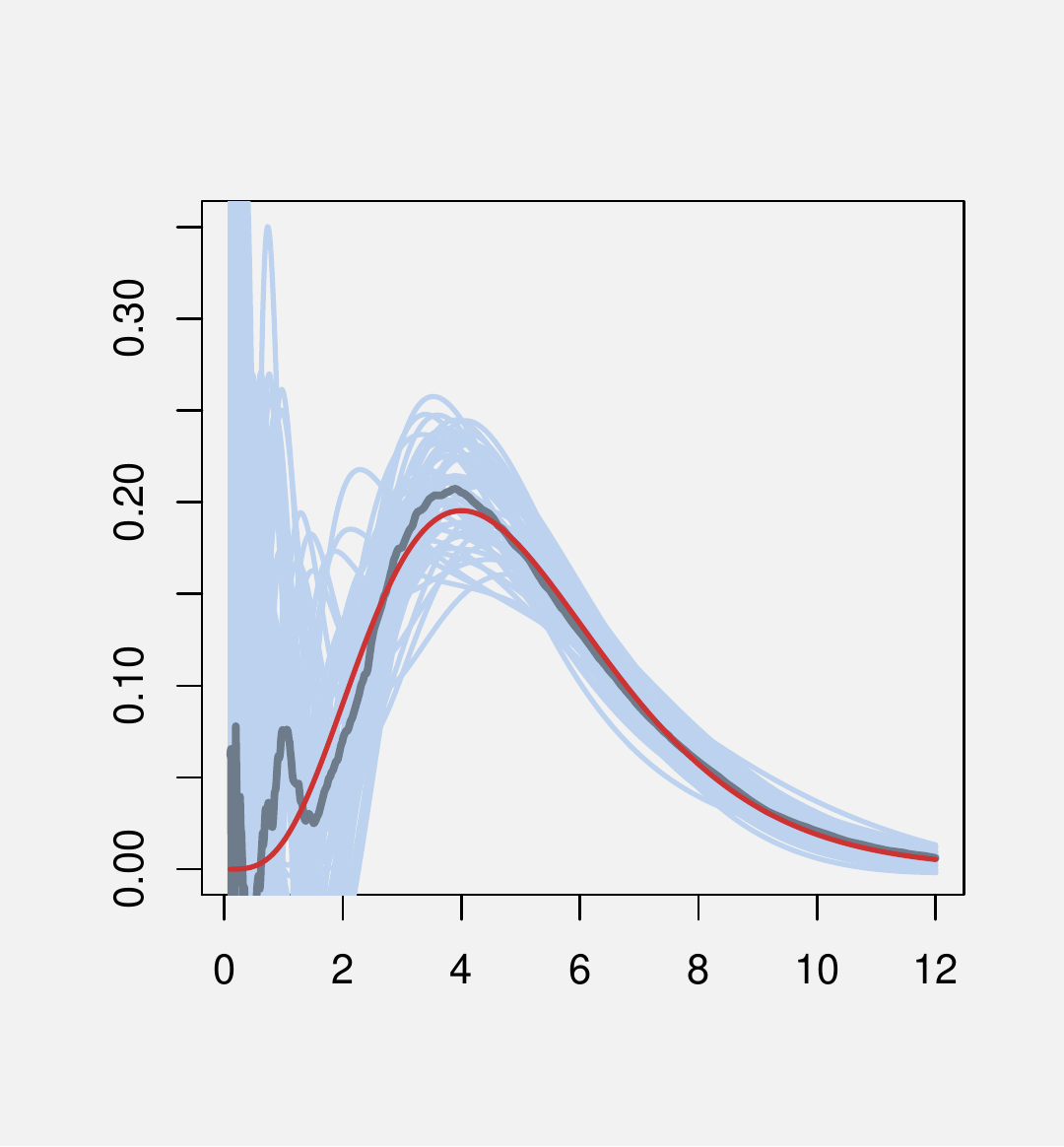}
	\end{minipage}}
\centerline{\begin{minipage}[t]{0.32\textwidth}
	\centering		\includegraphics[width=\textwidth,height=35mm]{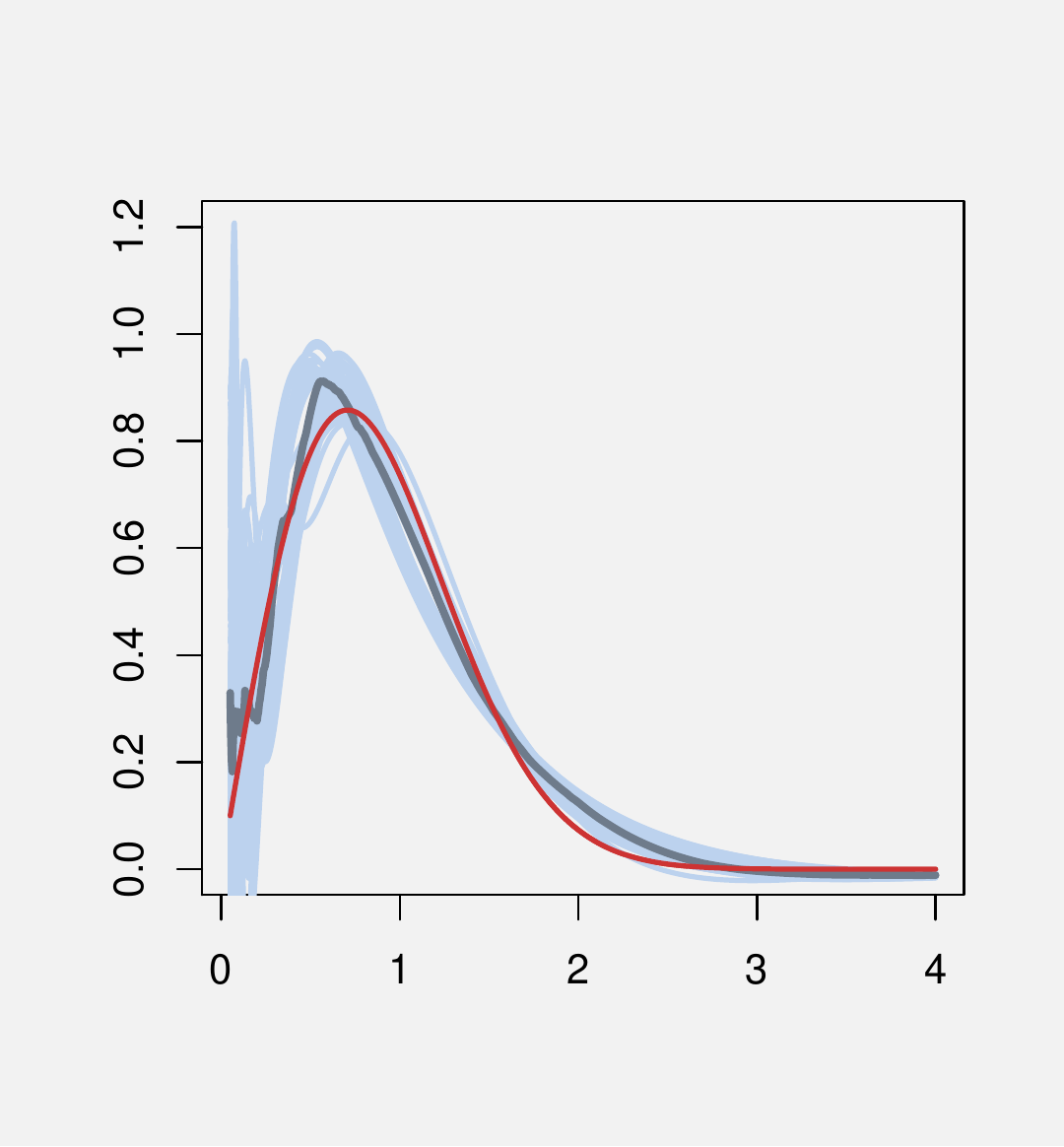}
	\end{minipage}
	\begin{minipage}[t]{0.32\textwidth}
	\centering		\includegraphics[width=\textwidth,height=35mm]{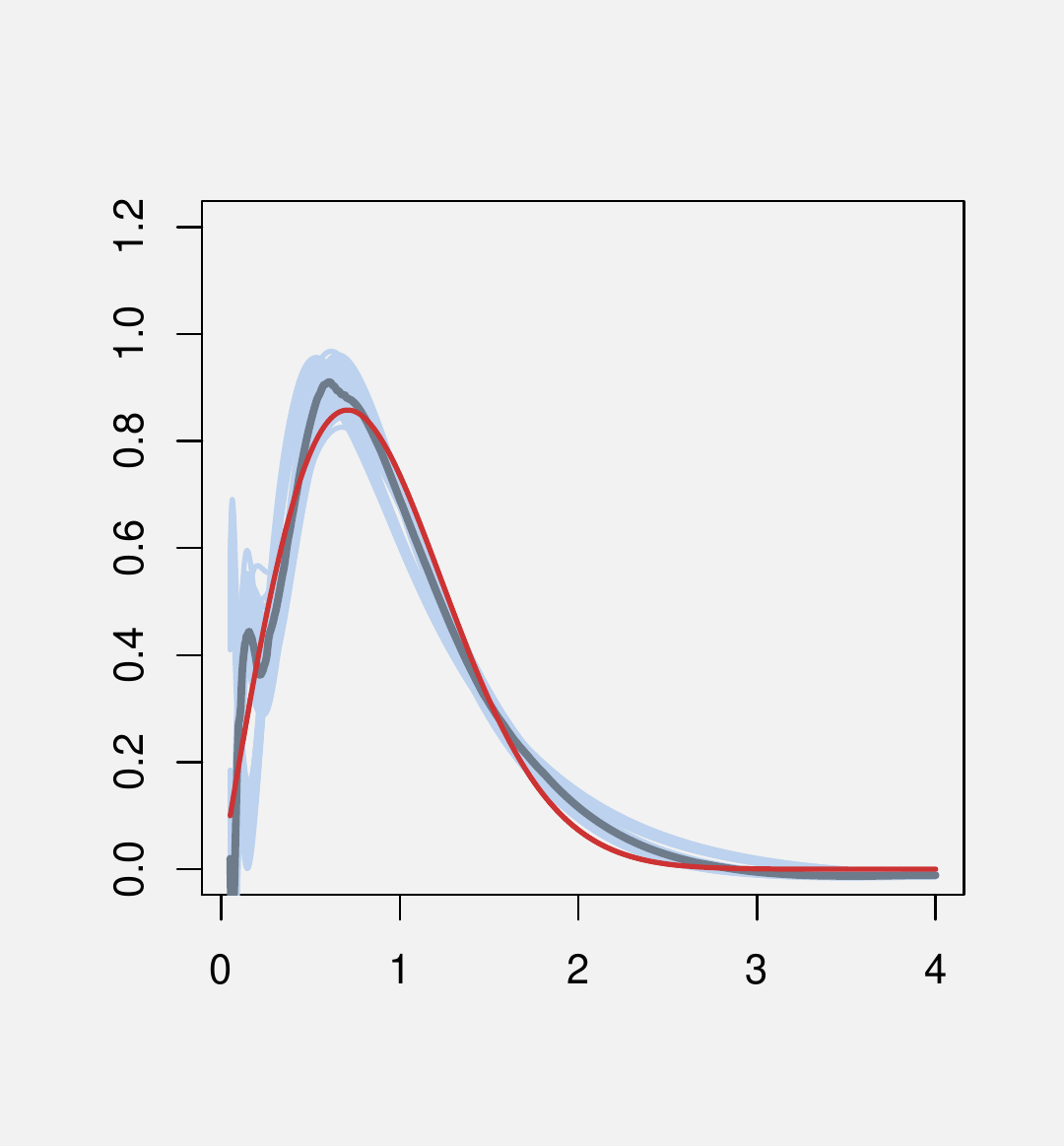}
	\end{minipage}
	\begin{minipage}[t]{0.32\textwidth}
	\centering		\includegraphics[width=\textwidth,height=35mm]{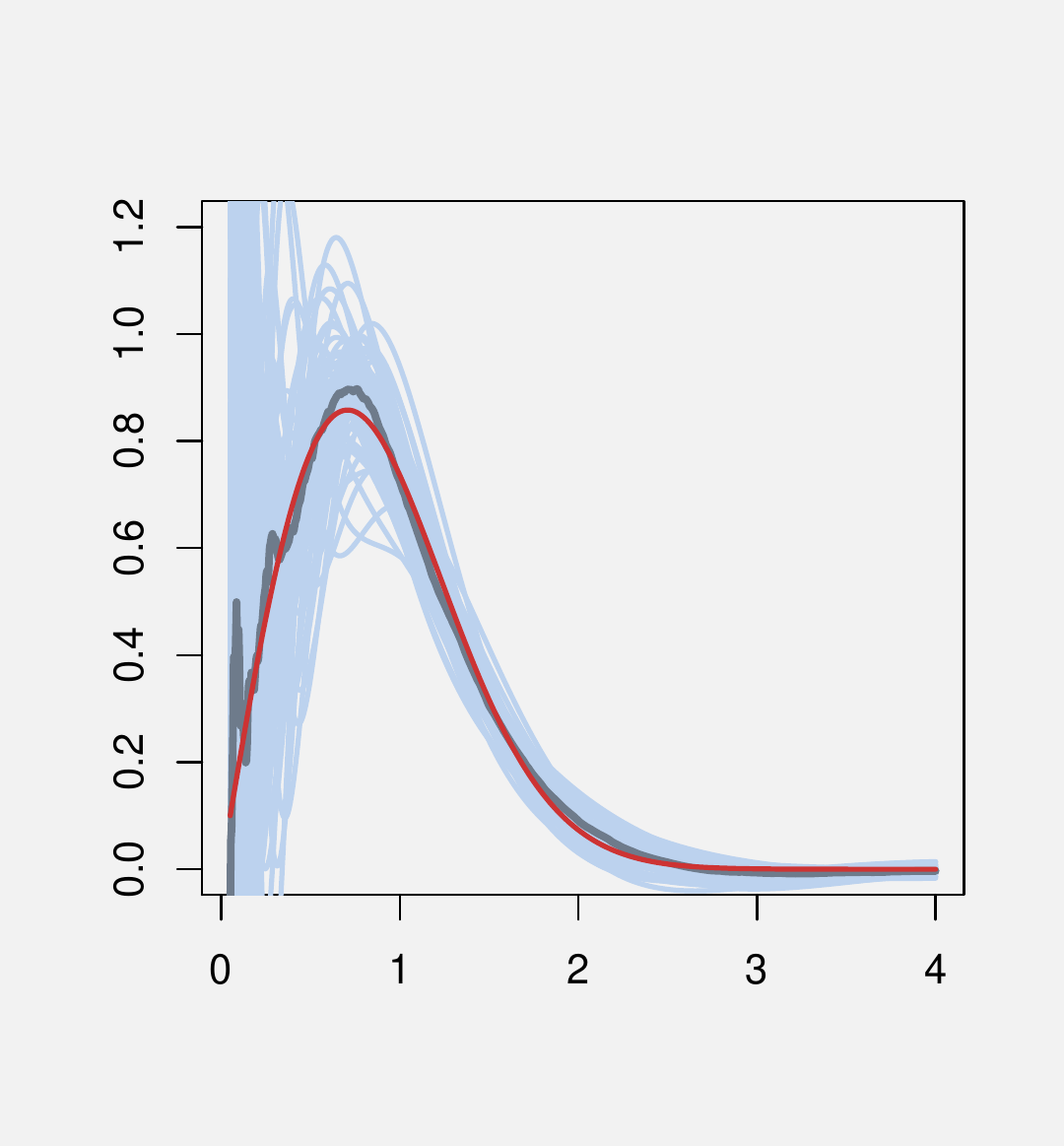}
	\end{minipage}}
\centerline{\begin{minipage}[t]{0.32\textwidth}
		\centering		\includegraphics[width=\textwidth,height=35mm]{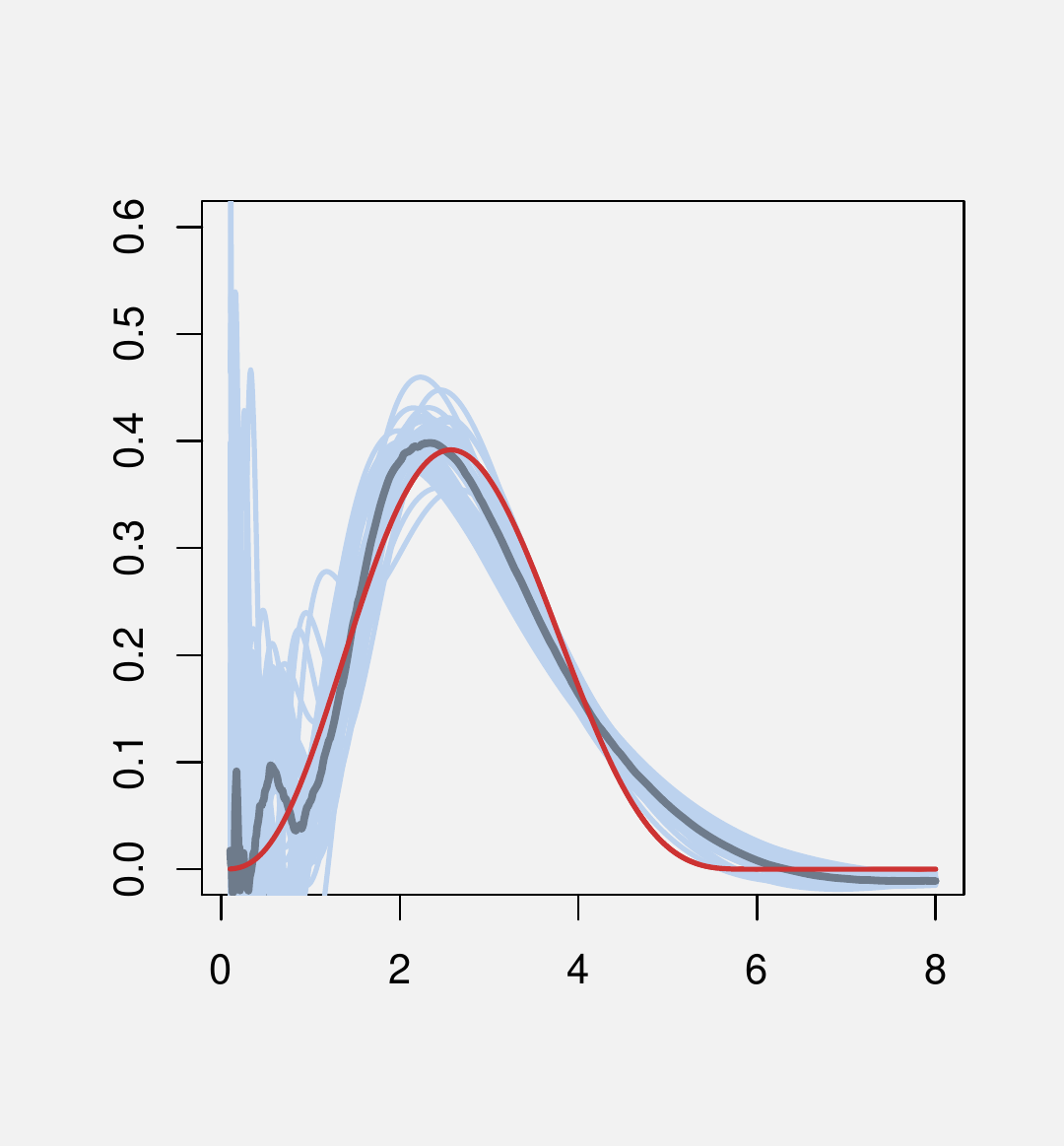}
	\end{minipage}
	\begin{minipage}[t]{0.32\textwidth}
		\centering		\includegraphics[width=\textwidth,height=35mm]{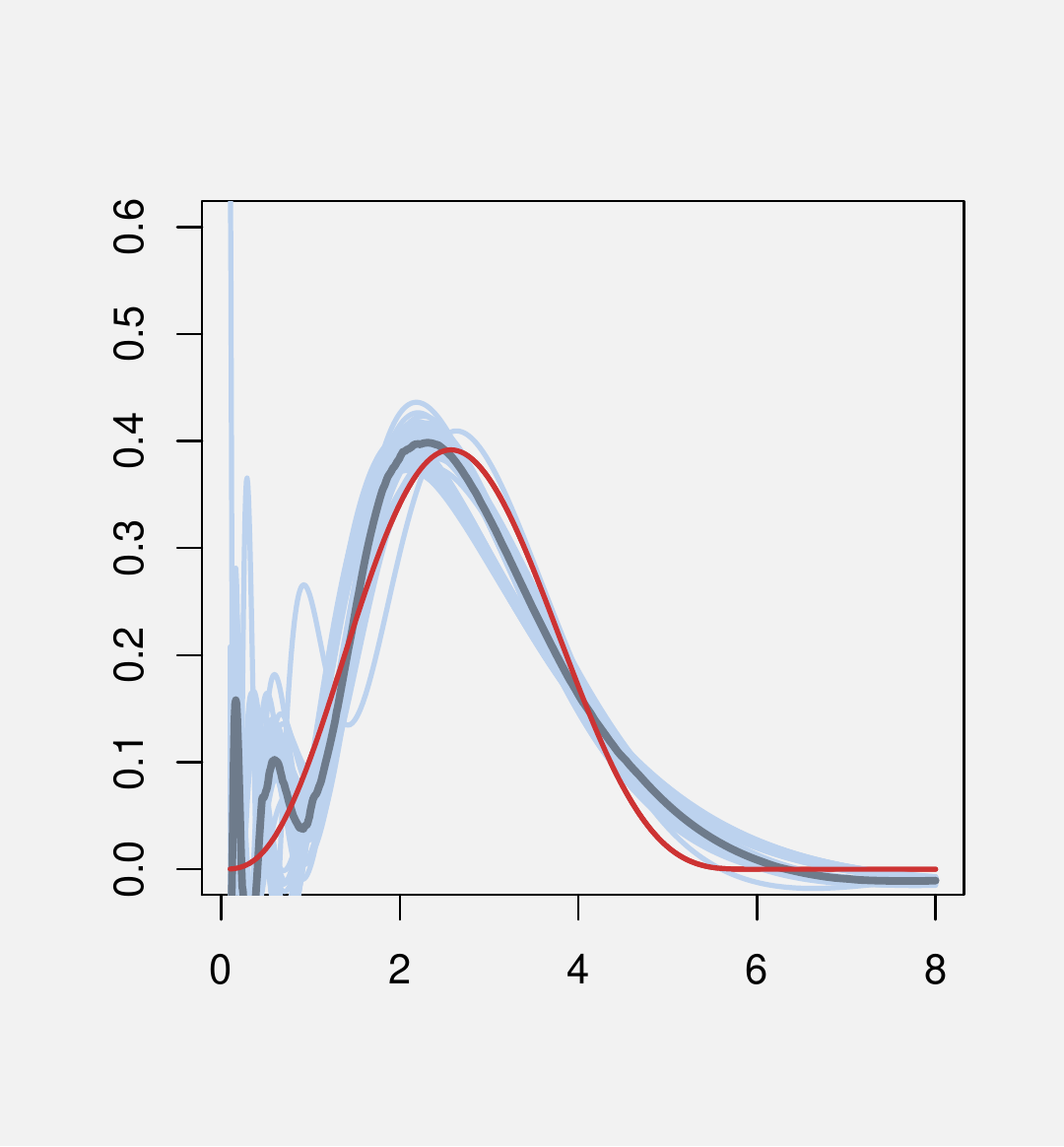}
	\end{minipage}
	\begin{minipage}[t]{0.32\textwidth}
		\centering		\includegraphics[width=\textwidth,height=40mm]{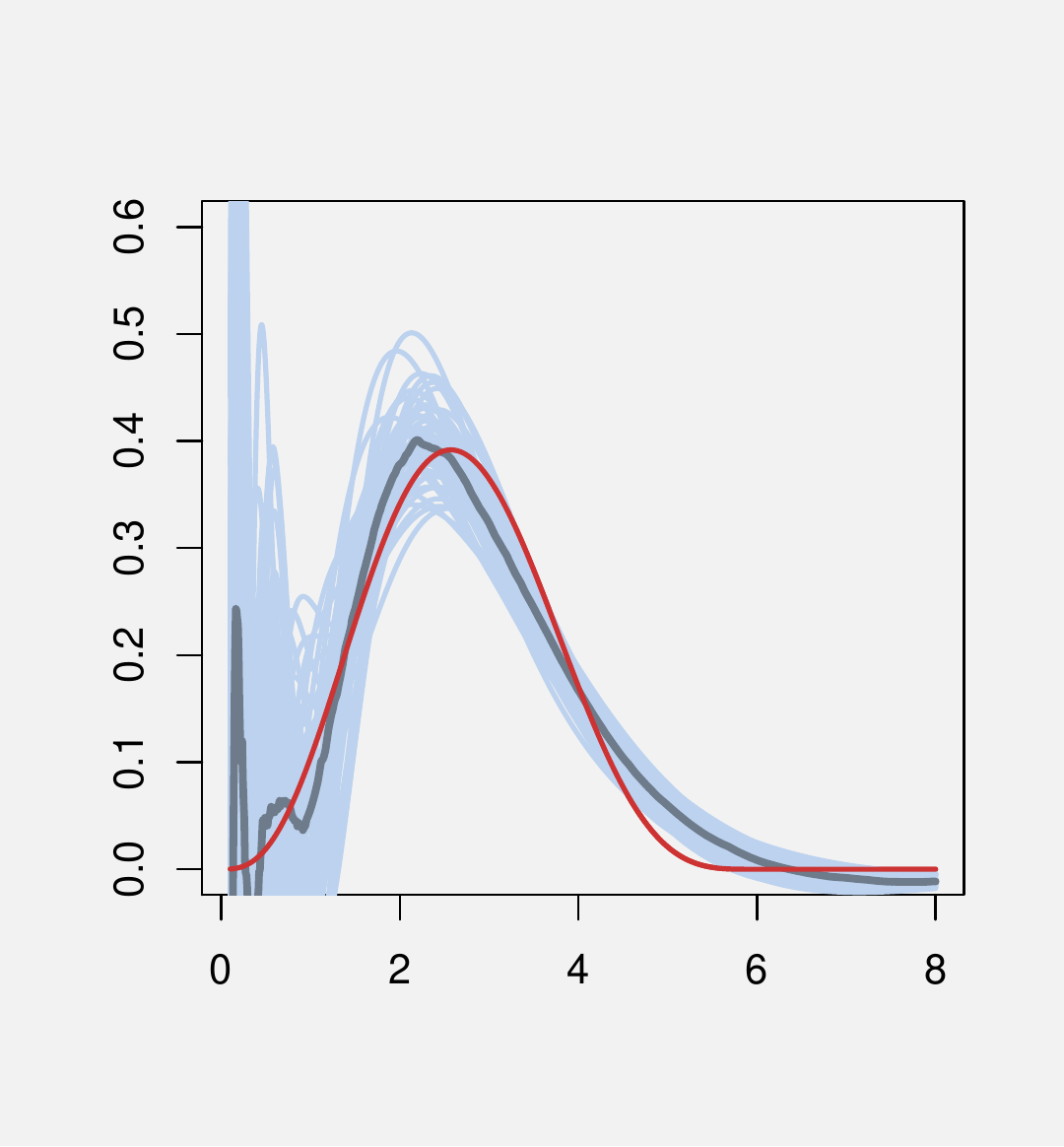}
\end{minipage}}
	\captionof{figure}{\label{figure:2}Considering the estimator $\widehat f_{\widehat k}$ and a sample size
          $n=2000$ the adaptive estimators are depicted for 50
          Monte-Carlo simulations with $U\sim U_{[0,s]}$ (left), $U\sim U_{[1/2, 3/2]}$ (middle) and $U\sim \beta(1,2)$ (right)
         in the cases \ref{si:lag:i} (first row), \ref{si:lag:iv} (second row) and \ref{si:lag:iii} (third row). The true density $\So$ is given by the red curve while the dark blue curve is the point-wise empirical median of the 50 estimates.}
\end{minipage}\\[2ex]
\paragraph{Comment}
Since the estimator $\widehat f_{\widehat k}$ is built to minimize the weighted global risk, it seems natural that the estimator $\widehat f_{\widehat k}$ behaves worse in the region close to zero then the estimator $\widetilde f$ which is built to minimize the unweighted global risk. This effect is observable in \cref{figure:1}. Furthermore, the developed minimax theory suggests that the cases of $U\sim \text{U}_{[0,1]}$ and $U\sim \text{U}_{[0.5,1.5]}$ are of similiar complexity which is reflected in the plots of \cref{figure:2}.
For the case $U\sim \beta(1,2)$, parameter $\gamma=2$ in \textbf{[G1]}, both theory and simulation imply that the recovering of the density $f$ based on the noise sample $Y_1,\dots,Y_n$ leads to a more difficult inverse problem than in the other cases.

%%% Local Variables:
%%% mode: latex
%%% TeX-master: "_0SCORDEMME"
%%% End:

% --------------------------------------------------------------------
% <<Appendix>>
% --------------------------------------------------------------------
\appendix
\setcounter{subsection}{0}
\section*{Appendix}
\numberwithin{equation}{subsection}  
\renewcommand{\thesubsection}{\Alph{subsection}}
\renewcommand{\theco}{\Alph{subsection}.\arabic{co}}
\numberwithin{co}{subsection}
\renewcommand{\thelem}{\Alph{subsection}.\arabic{lem}}
\numberwithin{lem}{subsection}
\renewcommand{\therem}{\Alph{subsection}.\arabic{rem}}
\numberwithin{rem}{subsection}
\renewcommand{\thepr}{\Alph{subsection}.\arabic{pr}}
\numberwithin{pr}{subsection}
%======================================================================================================================
%                                                                 
% Title: Appendix 1
% Author: Sergio Brenner Miguel and Jan JOHANNES, IAM, Ruprecht-Karls Universität Heidelberg, Deutschland  
% 
% Date: %%ts latex start%%[2020-01-29 Wed 15:44]%%ts latex end%%
%
% ======================================================================================================================
\subsection{Preliminaries}\label{a:prel}
%

% --------------------------------------------------------------------
% <<Re \ref{re:diff_fun}>>
% --------------------------------------------------------------------

\paragraph{Properties of the Mellin transform}
By assuming that $h\in \Lz_{\pRz}^{1, \text{loc}}$ is a at least $b$-time differentiable function $h$, where $h^{(b)}$ denotes its $b$-th derivative, $ b\in \Nz$, and that $c-b\in \Xi_h$ and $c+a\in \Xi_h, a\in \Nz,$ we get that 
\begin{align}\label{equation:Mell_rule}
\hspace*{-0.5cm}\Mela {x^a h}{c}(t)= \Mela h {c+a}(t) \, \text{ respec. } \, \Mela{h^{(b)}}{c}(t) = (-1)^b \frac{\Gamma(c+it)}{\Gamma(c-b+it)} \Mela{h}{c-b}(t)
\end{align}
Combining both results in \eqref{equation:Mell_rule} we get that $\Mela{xh^{(1)}}c(t)=(-c-it)\Mela {h}{c}(t)$ if  $c+1\in \Xi_h$ and $h$ differentiable. Further for $h_1,h_2\in \LpA$ with $c\in \Xi_{h_1}$ and $1-c\in \Xi_{h_2}$ we get that $\int_0^{\infty} h_1(x)h_2(x)dx = \frac{1}{2\pi} \int_{-\infty}^{\infty} \Mela {h_1} c(t) \overline{\Mela {h_2} {1-c}(t)} dt$. Combining this and \eqref{equation:Mell_rule} we conclude that for  $h_1, h_2 \in \mathbb L_2(\omega_{\alpha})$ it holds $\langle h_1, h_2 \rangle_{\omega_{\alpha}}= \frac{1}{2\pi}\int_{-\infty}^{\infty} \Mela{h_1}{\alpha}(t) \overline{\Mela{h_2}{\alpha}(t)}dt$.
\begin{lem}(Talagrand's inequality)\label{tal:re} Let
	$\Ob_1,\dotsc,\Ob_n$ be independent $\cZ$-valued random variables and let $\bar{\nu}_{h}=n^{-1}\sum_{i=1}^n\left[\nu_{h}(\Ob_i)-\Ex\left(\nu_{h}(\Ob_i)\right) \right]$ for $\nu_{h}$ belonging to a countable class $\{\nu_{h},h\in\cH\}$ of measurable functions. Then,
	\begin{align}
	&\Ex\vectp{\sup_{\He\in\cH}|\bar{\nu_{h}}|^2-6\TcH^2}\leq C \left[\frac{\Tcv}{n}\exp\left(\frac{-n\TcH^2}{6\Tcv}\right)+\frac{\Tch^2}{n^2}\exp\left(\frac{-K n \TcH}{\Tch}\right) \right]\label{tal:re1} 
	\end{align}
	with numerical constants $K=({\sqrt{2}-1})/({21\sqrt{2}})$ and $C>0$ and where
	\begin{equation*}
	\sup_{\He\in\cH}\sup_{z\in\cZ}|\nu_{h}(z)|\leq \Tch,\qquad \Ex(\sup_{\He\in\cH}|\bar{\nu_{h}}|)\leq \TcH,\qquad \sup_{\He\in\cH}\frac{1}{n}\sum_{i=1}^n \Var(\nu_{h}(\Ob_i))\leq \Tcv.
	\end{equation*}
\end{lem}
\begin{rem}
	Keeping the bound \eqref{tal:re1}  in mind, let us
	specify particular choices $K$, in fact $K\geq \tfrac{1}{100}$.
	The next bound is now an immediate consequence, 
	\begin{align}
	&\Ex\vectp{\sup_{\He\in\cH}|\bar{\nu_{h}}|^2-6\TcH^2}\leq C \left(\frac{\Tcv}{n}\exp\left(\frac{-n\TcH^2}{6\Tcv}\right)+\frac{\Tch^2}{n^2}\exp\left(\frac{-n \TcH}{100\Tch}\right) \right)\label{tal:re3} 
	\end{align}
	In the sequel we will make use of the slightly simplified bounds \eqref{tal:re3} rather than \eqref{tal:re1}.
	\remEnd
\end{rem}

%%% Local Variables:
%%% mode: latex
%%% TeX-master: "_SCORDEMME"
%%% End:

%======================================================================================================================
%                                                                 
% Title: Appendix 3
% Author: Sergio Brenner Miguel and Jan JOHANNES, IAM, Ruprecht-Karls Universität Heidelberg, Deutschland
% 
% Date: %%ts latex start%%[2020-07-30 Thu 19:49]%%ts latex end%%
%
% ======================================================================================================================
\subsection{Proofs of \cref{ag}}\label{a:ag}
% --------------------------------------------------------------------
% <<Proof of Re key argument>>
% --------------------------------------------------------------------

\begin{pro}[Proof of \cref{dd:pr:con}]
Since $\Mela {f-f_k}{\alpha}(t)=0$ for $|t|\geq k$ we get that $\langle f-f_k, f_k -\widehat f_k \rangle_{\omega_{\alpha}} = \frac{1}{2\pi} \int_{-k}^{k} \Mela {f-f_k}{\alpha}(t) \overline{\Mela {f_k-\widehat f_k}{\alpha}(t)} dt =0$ and thus $\|f-\widehat f_k\|_{\omega_{\alpha}}^2 = \|f- f_k\|_{\omega_{\alpha}}^2 + \|f_k-\widehat f_k\|_{\omega_{\alpha}}^2$. 

Finally we see $\|\widehat f_k-f_k \|_{\omega_{\alpha}}^2 = \frac{1}{2\pi} \int_{-k}^k |\Mela f {\alpha} (t)-\widehat \M_{\alpha}(t)|^2 dt$ and $\E_{f}^n(|\Mela f{\alpha}(t)-\widehat \M_{ \alpha}(t)|^2) = \Var(\widehat \M_{\alpha}(t)) \leq \frac{1}{n}  \E_{f}^{1}({|X_1^{\alpha-1+it}|^2}) = \sigma^2 n^{-1}$. Now by using Fubini we get
$
\E_{f}^n(\|\widehat f_k-f_k\|_{\omega_{\alpha}}^2) \leq \sigma^2\pi^{-1} kn^{-1}.
$ \proEnd
\end{pro}

\begin{pro}[Proof of \cref{mm:pr:reg}]
	The main strategy of this proofs relies on the well-known fact that 
		\begin{align*}
	W^s=\{H\in \Lz_{\Rz}^2: \text{H weakly differentiable up } \text{ to the order } s, H^{(i)} \in \Lz_{\Rz}^2, i\in \nset{0,s} \}
	\end{align*}
	for $s\in \Nz$ and the already discussed connection between the Mellin transform and the Fourier transform. 
	Further we want to stress out that for a function $H\in \Lz_{\Omega}^2, \Omega\subset \Rz$ open,  being weakly differentiable corresponds to being locally absolutely continuous on $\Omega$, that is $h$ is absolutely continuous on all compact intervals $[a,b] \subset \Omega$, $a<b\in\Omega$. \\
	Let us start by assuming $f\in \Wz^s$ then $F:=(\omega f)\circ \varphi \in W^s$ which means that $F$ is $s$-times weakly differentiable with $F^{(i)} \in \Lz_{\Rz}^2$ for $i\in \nset{0,s}$. Thus we get that $F$ is a $s-1$-times continously differentiable function, more precisely there exists a representant of the equivalence class of $F$ such that it is $s-1$-times continuously differentiable. Now since $f=\omega^{-1} F(\varphi^{-1})$ we can deduce that $f$ itself is $s-1$-times continuously differentiable. Let us define the operator $\mathcal T: C^{1}(\pRz) \rightarrow C^{0}(\pRz), f\mapsto (\omega f)^{(1)}$, $\mathcal T^{0}:= \mathrm{Id}$ denote the identity and $\mathcal T^{j}= \mathcal T \circ \mathcal T^{j-1}$ for $j\in \Nz$. Since $\mathcal T[\omega^{j}f^{(j)}]=(j+1) \omega^j f^{(j)}+ \omega^{j+1} f^{(j+1)}$ for $j\in \nsetro{0,s-1}$ we conclude $\mathcal T^j[f]= \sum_{i=0}^j b_{i,j}\,\omega^i f^{(i)}$, $b_{i,j} \geq 1$. Now we can use the following Lemma to deduce  that $\| \mathcal T^j[f]\|_{\omega}^2 = \| (\omega\mathcal T^j[f]) \circ \varphi \|_{\Rz}^2  =\| F^{(j)}\|_{\Rz}^2 <\infty.$
	\begin{lem}\label{lem:help1}
		For $h\in C^k(\pRz)$ and all $j\in \nset{0, k}$, it holds $(\omega\mathcal T^j[h]) \circ \varphi = (-1)^{j}(\varphi h(\varphi))^{(j)}$.
	\end{lem}
	From this follows directly that $\omega^{j}f^{(j)} \in \LpA(\omega)$ for $j\in\nsetro{0,s}$.
	As the next step we show that $f^{(s-1)}$ is locally absolutely continuous. To do so, we see first for $j\in\nsetro{0,s}$ we have that $\mathcal T^j[f] \in \Wz^1$. Now using the following lemma implies that $\omega^{s-1} f^{(s-1)}$ is absolutely continuous as linear combination of $\{\mathcal T^j[f]:j\in\nsetro{0,s}\}$.
	\begin{lem}\label{lem:help2}
		Let $h\in \Wz^1$. Then $h$ is a locally absolutely continuous function with derivative $h':\pRz \rightarrow \Rz$ and $h, \omega h' \in \LpA(\omega)$.
	\end{lem}
	Let now $\delta$ denote the derivative of $\omega^{s-1}f^{(s-1)}$. Then by \cref{lem:help2} we see that $\omega \delta \in \LpA(\omega)$. Defining now $f^{(s)}:=\omega^{-s+1}\delta  -(s-1)\omega^{-1}f^{(s-1)} \in \Lz_{\text{loc}}^1(\pRz)$ we get for any $a,b\in \pRz$, $a<b$,
	$
	\int_a^b f^{(s)}(x)dx = f^{(s-1)}(b)-f^{(s-1)}(a)
	$
	using the integration by part rule for absolutely continuous function (see \cite{Cohn2013}). Finally we have that $\omega^sf^{(s)}=\delta \omega -(s-1) \omega^{s-1}f^{(s-1)} \in \LpA(\omega).$ \\
	Let us now show the other direction, indeed let us assume that $f$ is $s-1$-times continuously differentiable, $f^{(s-1)}$ is locally absolutely continuous with derivative $f^{(s)}$ and $\omega^{j} f^{(j)}\in \LpA(\omega)$ for $j\in\nset{0,s}$.
	Thus for  $\mathcal T^{j}[f]= \sum_{i=0}^j c_{i,j} \omega^i f^{(i)} \in\LpA(\omega)$ with $j\in\nsetro{0,s-1}$ we have that $\omega(\mathcal T^j[f])^{(1)}= \mathcal T^{j+1}[f]- \mathcal T^{j}[f]\in \LpA(\omega)$. We can conclude that $\mathcal T^{j}[f]\in \Wz^1$ for $j\in\nsetro{0,s-1}$ applying the following lemma.
	\begin{lem}\label{lem:help3}
		Let $h: \pRz \rightarrow \Rz$ be locally absolutely continous function with derivative $h': \pRz\rightarrow \Rz$ and $h, \omega h' \in \LpA(\omega)$. Then $h\in \Wz^1$.
	\end{lem}
	Now setting $\delta = (s-1)\omega^{s-2} f^{(s-1)}+ \omega^{s-1}f^{(s)} \in \Lz_{\text{loc}}^1(\pRz)$ we have that $\int_a^b \delta(x)dx = b^{s-1} f^{(s-1)}(b)- a^{s-1} f^{(s-1)}(a)$ and $\omega\delta\in \LpA(\omega)$. Thus again applying \cref{lem:help3}  on $\omega^s f^{(s)}$ again shows that $\omega^s f^{(s)} \in \Wz^1$ and thus $\mathcal T^{s-1}[f]\in \Wz^1$. Now we use that $\Mela{\mathcal T^{s-1}}{1}(t)=(-1)^{s-1}\mathcal F[F^{(s-1)}](t)=(-it)^{s-1} \mathcal F[F](t)= (-it)^{s-1}\Mela {f}{1}(t)$ which implies that $f\in \Wz^s$. 
	\proEnd
\end{pro}

\begin{pro}[Proof of \cref{lem:help1}]
	For $j=0$ the claim is trivially correct. Assume that the claim hold for $j\in \nsetro{1,k}$ then
	\begin{align*}
	(\omega\mathcal T^{j+1}[f]) \circ \varphi= (\omega\mathcal T^{j}[f+\omega f^{(1)}])\circ \varphi 
	= (-1)^{j-1} (\varphi f(\varphi) +\varphi^2 f^{(1)}(\varphi))^{(j)}
	\end{align*}
	and thus $(\varphi f(\varphi))^{(1)}= - \varphi f(\varphi)- \varphi^2 f^{(1)}(\varphi)$ implies the claim.
\end{pro}

\begin{pro}[Proof of \cref{lem:help2}]
	Since $h\in \Wz^1$ we have that $H=(\omega h) \circ \varphi$ lies in the Sobolev space of order 1. In equal $H$ is locally absolutely continuous with derivative $H'$ and $H, H' \in \mathbb L^2_{\Rz}$. 
	From this we can conclude that $h$ is locally absolutely continuous. Indeed for $h':= -\omega^{-2} (H'+H)\circ \varphi^{-1} \in \mathbb L_{\text{loc}}^1(\pRz)$ and $a,b\in\pRz, a<b$ holds
	\begin{align*}
	\hspace*{-1cm}\int_a^{b} -x^{-2} H'(\varphi^{-1})(x)dx = \int_{\varphi^{-1}(a)}^{\varphi^{-1}(b)} \exp(x) H'(x) dx = h(b)-h(a) + \int_a^b x^{-2} H(\varphi^{-1})(x) dx
	\end{align*}
	applying the integration by part rule for absolutely continuous function. Further we have that $\|\omega h'\|_{\omega}= \|\omega^{-1}(H'+H) \circ \varphi^{-1}\|_{\omega} \leq \| H'\|_{\Rz}+\| H\|_{\Rz}< \infty.$ 
\end{pro}

\begin{pro}[Proof of \cref{lem:help3}]
	Since $h\in \LpA(\omega)$ we have for $H:= (\omega h)\circ \varphi$ that $\|H\|_{\Rz}=\|h\|_{\omega} <\infty.$ Further $H$ is locally absolutely continuous with derivative $-\varphi h(\varphi)-\varphi^2 h'(\varphi)$ since for $a,b\in \Rz$ with $a<b$ holds
	\begin{align*}
	\int_a^{b} -\varphi^2(x) h'(\varphi(x))dx &= \int_{\varphi(a)}^{\varphi(b)} x h'(x) dx 
	= \left[\varphi h \circ \varphi \right]_a^b - \int_{\varphi(a)}^{\varphi(b)} h(x)dx. 
	\end{align*}
	Now since $\|\varphi^2 h'(\varphi)\|_{\Rz}=\|\omega f'\|_{\omega}<\infty$ we deduce that $H$ is in the Sobolev space of order 1 and thus $(1+t^2)^{1/2}\Mela h 1 = (1+t^2)^{1/2}\mathcal F[H] \in \mathbb L^2(\Cz)$ and thus $h\in \Wz^1$.
\end{pro}

\begin{pro}[Proof of \cref{dd:thm:ada}] 
	Let us define the nested subspaces $(S_k)_{k\in \pRz}$ by $S_k:=\{h\in \LpA(\omega): \forall |t|\geq k: \Mel[h](t)=0\}$. For any $h \in S_k$ we consider the empirical contrast 
	\begin{align*}
	\gamma_n(h) = \|h\|_{\omega}^2 -2  \frac{1}{2\pi} \int_{-\infty}^{\infty} \widehat{\M}(t)\frac{\Mel[h](-t)}{\Mel[g](t)}dt= \|h\|_{\omega}^2 - 2 n^{-1} \sum_{j=1}^n \nu_h(Y_j) 
	\end{align*}
	 with $\nu_h(Y_j):= \frac{1}{2\pi} \int_{-\infty}^{\infty} Y_j^{it} \frac{\Mel[h](-t)}{\Mel[g](t)}dt$. One can easily see that $\widehat f_k = \argmin\{\gamma_n(h):h\in S_k\}$ with $\gamma_n(\widehat f_k)=-\|\widehat f_k\|_{\omega}^2$.
	For $h\in S_k$  define the  empirical process $\bar{\nu}_h:= n^{-1} \sum_{j=1}^n \nu_h(Y_j) - \langle h, f \rangle_{\omega}$. Then we have that for $h_1,h_2 \in S_k$ that
	\begin{align}\label{dd:eq:cont}
	\gamma_n(h_1)-\gamma_n(h_2) = \|h_1-f\|_{\omega}^2 - \|h_2-f\|_{\omega}^2 - 2 \bar{\nu}_{h_1-h_2}.
	\end{align}
	Now since $\gamma_n(\widehat f_k) \leq \gamma_n(f_k)$ we get $\|f-\widehat f_k \|_{\omega}^2 \leq \|f-f_k\|^2_{\omega} + 2 \bar\nu_{\widehat f_k -f_k}$. 
	By definition of $\widehat k$ we have that $\gamma_n(\widehat f_{\widehat k})- \mathrm{pen}(\widehat k) \leq \gamma_n(\widehat f_{ k})- \mathrm{pen}( k) \leq \gamma_n(f_k)- \mathrm{pen}(k)$ for any $k\in \mathcal K_n$. Now using \eqref{dd:eq:cont} we get that 
	\begin{align*}
	\| f- \widehat f_{\widehat k} \|_{\omega}^2 \leq \|f-f_k\|_{\omega}^2 + 2 \bar\nu_{\widehat f_{\widehat k}- f_k} + \mathrm{pen}( k)- \mathrm{pen}(\widehat k).
	\end{align*}
	First we note that $S_{k_1} \subseteq S_{k_2}$ for $k_1 \leq k_2$. Let us now denote by $a\vee b:= \max(a,b)$ and define for all $k\in \mathcal K_n$ the unit balls $B_k:=\{h\in S_k: \|h\|_{\omega} \leq 1\}$. Next we deduce from $2ab \leq a^2+b^2$ that $ 2\bar\nu_{\widehat f_{\widehat k}-f_k}\leq 4^{-1}\|\widehat  f_{\widehat k} - f_k\|_{\omega}^2 + 4 \sup_{h \in B_{\widehat k \vee k} }\bar \nu_h^2$. Further we see that $4^{-1} \|\widehat  f_{\widehat k} - f_k\|_{\omega}^2  \leq 2^{-1} (\|\widehat  f_{\widehat k} - f\|_{\omega}^2 + \|f- f_k\|_{\omega}^2)$. Putting all this together and define 
	\begin{align}\label{eq:delta}
	p(\widehat k \vee k):= 6(2\pi n)^{-1}\Delta_g(\widehat k \vee k) \text{ where }\Delta_g(k):= \int_{-k}^k |\Mel[g](t)|^{-2} dt
	\end{align} we get
	\begin{align*}
	\hspace*{-1cm}	\| f- \widehat f_{\widehat k} \|_{\omega}^2 \leq 3\|f-f_k\|_{\omega}^2 + 8 \big(\sup_{h \in B_{\widehat k \vee k} }\bar\nu_h^2 -p(k\vee \widehat k) \big)_+ +8p(\widehat k \vee k) + 2\mathrm{pen}(k)- 2\mathrm{pen}(\widehat k)
	\end{align*}
	Assuming now that $\chi \geq 12 C_g \pi^{-1}$ we get that $4p(\widehat k \vee k) \leq \mathrm{pen}(k)+ \mathrm{pen}(\widehat k)$ and thus
	\begin{align*}
	\| f- \widehat f_{\widehat k} \|_{\omega}^2 \leq 4\big(\|f-f_k\|_{\omega}^2 +\mathrm{pen}(k) \big) + 8 \max_{k'\in \mathcal K_n}\big(\sup_{h \in B_{ k'} }\bar\nu_h^2 -p(k') \big)_+ 
	\end{align*}
	We will use the following lemma which we will be proven afterwards.
		\begin{lem}\label{dd:lem:talapply}
		Assuming that $\|\omega f_Y\|_{\infty}< \infty$ and that for all $k\in \mathcal K_n$ the function $G_k:\Rz \rightarrow \Rz, t\mapsto \1_{[-k,k]}(t) |\Mel[g](t)|^{-2}$ is bounded  we have 
		\begin{align*}
		\hspace*{-1cm}\E_{f_Y}^n\big(\sup_{h \in B_{ k} }\bar \nu_h^2 -p(k) \big)_+  \leq \frac{C}{n} &\left( \|G_k\|_{\infty}\|\omega f_Y\|_{\infty} \exp(-\frac{\Delta_g(k)}{12 \pi\|\omega f_Y\|_{\infty} \|G_k\|_{\infty} }) \right.\\
		&\left. +  \frac{\Delta_g(k)}{(2\pi)^2n} \exp(-\frac{\sqrt{n}}{50})  \right),
		\end{align*}
		where $\Delta_g$ is defined in \eqref{eq:delta}.
	\end{lem}
	Now under \textbf{[G1]} we have that $\Delta_g(k) \geq c_g k^{2\gamma+1} $ and for all $t\in \Rz$ holds $|G_k(t)| \leq C_g k^{2\gamma}$ thus we have that the first summand is bounded by $C_g k^{2\gamma}\|\omega f_Y\|_{\infty} \exp(-\frac{c_g k}{12 \pi\|\omega f_Y\|_{\infty} C_g })$ which is bounded over $\Nz$. For the second summand we use that $n^{-1} \Delta_g(k) \leq C_g n^{-1}k^{2\gamma+1} \leq C_g$ and thus bounded in $\Nz$.
	Applying the lemma we get that
	\begin{align*}
	\E_{f_Y}^n (	\| f- \widehat f_{\widehat k} \|_{\omega}^2) \leq 4\big(\|f-f_k\|_{\omega}^2 +\mathrm{pen}(k) \big) + C(\|\omega f_Y\|_{\infty},g) n^{-1}.
	\end{align*}
	Since this inequality holds for all $k\in \mathcal K_n$ this implies the claim.
	\proEnd
\end{pro}
	 
\begin{pro}[Proof of \cref{dd:lem:talapply}]
	We will use the Talagrand inequality \eqref{tal:re3} to show the claim. We want to emphasize that we are able to apply the Talagrand inequality on the sets $B_k$ since $B_k$ has a dense countable subset and due to continuity arguments. To do so we start to determine the constant $\Psi^2$. We have for any $h\in B_{k}$ that $\bar \nu_h^2 = \langle h, \widehat f_{k} - f_{k} \rangle_{\omega}^2\leq \| h\|_{\omega}^2 \|\widehat f_{k}- f_{k}\|_{\omega}^2$. Since $\|h\|_{\omega}\leq1$ we get 
	\begin{align*}
	\E_{f_Y}^n( \sup_{h\in B_{k}} \bar\nu_h^2 ) \leq \E_{f_Y}^n(\|\widehat f_{k} - f_{k}\|_{\omega}^2) \leq (2n \pi)^{-1}\Delta_g(k)=:\Psi^2.
	\end{align*}
	Thus $6\Psi^2 = p(k)$.
	Next we consider $\psi$. Let $y>0$ and $h\in B_{k}$ then using the Cauchy Schwartz inequality we get $|\nu_h(y)|^2= (2\pi)^{-2}|\int_{-k}^{k} y^{it} \frac{\Mel[h](-t)}{\Mel[g](t)} dt |^2 \leq (2\pi)^{-2} \int_{-k}^{k}  |\Mel[g](t)|^{-2}dt\leq (2\pi)^{-2} \Delta_g(k)=: \psi^2$ since $|y^{it}|=1$ for all $t\in \Rz$.\\
	Next we consider $\tau$. In fact for $h\in B_{k}$ we can conclude $\Var(\nu_h(Y_1)) \leq \E_{f_Y}^n(\nu_h(Y_1)^2)  \leq \|\omega f_Y\|_{\infty}
	\int_0^{\infty} y^{-1} \nu_h(y)^2 dt= \|\omega f_Y\|_{\infty}\| \nu_h\|_{\omega_{0}}^2$ with $\nu_h(y)= (2\pi)^{-1}\int_{-k}^k y^{it}\frac{\Mel[h](-t)}{\Mel[g](t)}dt$ for $y>0$. Thus 
	\begin{align*}
	\| \nu_h\|_{\omega_{0}}^2=\frac{1}{2\pi} \int_{-k}^{k} \left| \frac{\Mel[h](t)}{\Mel[g](t)}\right|^2 dt \leq \frac{\|G_k\|_{\infty}}{2\pi} \int_{-\infty}^{\infty} |\Mel[h](t)|^2 dt
	\end{align*}
	where $\frac{1}{2\pi} \int_{-\infty}^{\infty} |\Mel[h](t)|^2 dt = \|h\|_{\omega}^2 \leq 1$. Thus we set $\tau = \|\omega f_Y\|_{\infty}\|G_k\|_{\infty}$.
	Hence we have that $\frac{n\Psi^2}{6\tau}= \frac{\Delta_g(k)}{12 \pi\|\omega f_Y\|_{\infty} \|G_k\|_{\infty} }$ and  $\frac{n\Psi}{\psi}=\sqrt{2\pi n}$. We deduce
	\begin{align*}
	\hspace*{-1cm}\E_{f_Y}^n\big(\sup_{h \in B_{ k} }\bar \nu_h^2 -p(k) \big)_+  \leq \frac{C}{n} &\left( \|G_k\|_{\infty}\|\omega f_Y\|_{\infty} \exp(-\frac{\Delta_g(k)}{12 \pi\|\omega f_Y\|_{\infty} \|G_k\|_{\infty} }) \right.\\
	&\left. +  \frac{\Delta_g(k)}{(2\pi)^2n} \exp(-\frac{\sqrt{n}}{50})  \right).
	\end{align*}
	
\end{pro}

%%% Local Variables:
%%% mode: latex
%%% TeX-master: "_0SCORDEMME"
%%% End:
 
%======================================================================================================================
%                                                                 
% Title: Appendix 2
% Author: Sergio Brenner Miguel and Jan JOHANNES, IAM, Ruprecht-Karls Universität Heidelberg, Deutschland
% 
% Date: %%ts latex start%%[2020-01-28 Tue 15:38]%%ts latex end%%
%
% ======================================================================================================================
\subsection{Proofs of \cref{mm}}\label{a:mm}

\begin{pro}[Proof of \cref{theorem:lower_bound}]
	First we outline here the main steps of the proof.  We will construct a family of
	functions in $\rwcSobD{\wSob,\rSob}$ by a perturbation of the
	density $\SoPr[o]: \pRz \rightarrow \pRz$ with small bumps, such that their $\LpA(\omega)$-distance
	and the Kullback-Leibler divergence of their induced distributions can be bounded from below and
	above, respectively. The claim follows then by applying Theorem 2.5
	in \cite{Tsybakov2008}. We use
	the following construction, which we present first.\\ 
	Denote by $C_c^{\infty}(\Rz)$  the set of all smooth functions with compact support in $\Rz$ and let $\psi\in C_c^{\infty}(\Rz)$ be a function with support in $[0,1]$ and $\int_0^{1} \psi(x)dx = 0$. For each $K\in\Nz$ (to be selected below) and
	$k\in\nsetro{0,K}$ we define the bump-functions
	$\psi_{k, K}(x):= \psi(xK-K-k),$ $x\in\Rz$. and define for $j\in \Nz_0$ the finite constant $C_{j,\infty}:= \max(\|\psi^{(l)}\|_{\infty}, l\in \nset{0,j})$. Let us further define the operator $\mathcal S: C_c^{\infty}(\Rz)\rightarrow C_c^{\infty}(\Rz)$ with $\mathcal S[f](x)=x f^{(1)}(x)$ for all $x\in \Rz$ and define $\mathcal S^1:=\mathcal S$ and $\mathcal S^{n}:=\mathcal S \circ \mathcal S^{n-1}$ for $n\in \Nz, n\geq 2$.   Now, for $j \in \Nz$, we define the function $
	\psi_{k,K,j}(x):= \mathcal S^{j} [\psi_{k,K}](x)=\sum_{i=1}^{j} c_{i,j} x^i K^{i} \psi^{(i)}(xK-K-k)$ for $x \in \pRz$ and $c_{i,j} \geq 1$ and let $c_j:= \sum_{i=1}^j c_{i,j}$ \\
	For a bump-amplitude $\delta>0, \gamma \in \Nz$ and a vector 
	$\bm{\theta}=(\theta_1,\dots,\theta_K)\in \{0,1\}^K$ we define
	\begin{equation}\label{equation:lobodens}
	\SoPr[\bm{\theta}](x)=\SoPr[o](x)+ \delta K^{-s-\gamma} \sum_{k=0}^{K-1}
	\theta_{k+1} \psi_{k, K,\gamma}(x)\text{ where } f_o(x):=\exp(-x).
	\end{equation}
	Until now, we did not give a sufficient condition to ensure that our constructed functions $\{f_{\bm{\theta}}: \bm{\theta} \in \{0,1\}^K\}$ are in fact densities. This condition is given by the following lemma.
	\begin{lem}\label{mm:lem:den}
		Let $0<\delta< \delta_o(\psi, \gamma):=\exp(-2)2^{-\gamma} (C_{\gamma,\infty} c_{\gamma})^{-1}$. Then for all $\bm{\theta}\in\{0,1\}^K$, $f_{\bm{\theta}}$ is a density.
	\end{lem}
	Further one can show that these densities  all lie inside the ellipsoids $\rwcSobD[\Rz^+]{\wSob,L}$ for $L$ big enough. This is captured in the following lemma.
		\begin{lem}\label{lemma:Lag_SobDen}Let
		$\wSob\in\Nz$. Then,
		there is $L_{\wSob,\gamma,\delta}>0$ such that $\SoPr[o]$
		and any $\SoPr[\bm{\theta}]$ as in  \eqref{equation:lobodens} with
		$\bm\theta\in\{0,1\}^K$, $K\in\Nz$,  belong to $\rwcSobD[\Rz^+]{\wSob,L_{\wSob,\gamma,\delta}}$.
	\end{lem}
	For sake of simplicity we denote for a function $\varphi \in \LpA$ the multiplicative convolution with $g$ by $\widetilde{\varphi} := [\varphi *g]$. Futher we see that for  $y_2 \geq y_1 >0$ holds
	\begin{align*}
	\widetilde{f_o}(y_1) = \int_0^{\infty} g(x) x^{-1} \exp(-y_1/x) dx \geq \int_0^{\infty} g(x) x^{-1} \exp(-y_2/x) dx = \widetilde{f_o}(y_2)
	\end{align*} and thus $\widetilde{f_o}$ is monotone decreasing. Further we have that $\widetilde f_o(2)>0$ since otherwise $g=0$ almost everywhere.
	Exploiting \textit{Varshamov-Gilbert's
	lemma} (see \cite{Tsybakov2008}) in \cref{lemma:tsyb_vorb} we show
	further that there is $M\in\Nz$ with $M\geq 2^{K/8}$ and a subset
	$\{\bm \theta^{(0)}, \dots, \bm \theta^{(M)}\}$ of $\{0,1\}^K$ with
	$\bm \theta^{(0)}=(0, \dots, 0)$ such that for all
	$j, l \in \nset{0, M}$, $j \neq l$ the $\LpA(\omega)$-distance and the
	Kullback-Leibler divergence are bounded for $K\geq K_o(\gamma,\psi)$.
	\begin{lem}\label{lemma:tsyb_vorb}
		Let $K\geq K_o(\psi,\gamma)\vee 8$. Then there exists a subset $\{\bm \theta^{(0)}, \dots, \bm \theta^{(M)}\}$ of $\{0,1\}^K$  with $\bm \theta^{(0)}=(0, \dots, 0)$ such that $M\geq 2^{K/8}$ and for all $j, l \in \llbracket 0, M \rrbracket, j \neq l$ holds $\| f_{\bm\theta^{(j)}}-f_{\bm\theta^{(l)}}\|^2_{\omega} \geq \frac{\|\psi^{(\gamma)}\|^2\delta^2}{16}  K^{-2s}$ and $ \text{KL}(\widetilde{f}_{\bm\theta^{(j)}}, \widetilde{f}_{\bm\theta^{(0)}}) \leq \frac{C_1(g)\|\psi\|^2}{\widetilde{f}_o(2)\log(2)} \delta^2 \log(M) K^{-2s-2\gamma-1}$ where $\text{KL}$ is the Kullback-Leibler-divergence.
	\end{lem}
	Selecting $K=\ceil{n^{1/(2s+2\gamma+1)}}$, it follows
	\begin{align*}
	\frac{1}{M}\sum_{j=1}^M
	\text{KL}((\widetilde{f}_{\bm{\theta^{(j)}}})^{\otimes
		n},(\widetilde{f}_{\bm{\theta^{(0)}}})^{\otimes n})
	&= \frac{n}{M} \sum_{j=1}^M \text{KL}(
	\widetilde{f}_{\bm{\theta^{(j)}}},\widetilde{f}_{\bm{\theta^{(0)}}} )
	\leq \cst[(2)]{\psi, \delta,g,\gamma,f_o}  \log(M)
	\end{align*}
	where $\cst[(2)]{\psi, \delta,g,\gamma,f_o}< 1/8$ for all
	if $\delta\leq \delta_1(\psi, g, \gamma, f_o)$ and $M\geq 2$ for
	$n\geq n_{s,\gamma}:=8^{2s+1}\vee K_o(\gamma, \psi)^{2s+2\gamma+1}$. Thereby, we can use Theorem 2.5 of
	\cite{Tsybakov2008}, which in turn for any estimator $\hSo$ of $\So$
	implies
	\begin{multline*}
	\sup_{\So\in\rwcSobD{\wSob,\rSob}}
	\nVg\big(\|\hSo-\So\|_{\omega}^2\geq
	\tfrac{\cst[(1)]{\psi, \delta, \gamma}}{2}n^{-2s/(2s+2\gamma+1)} \big)\geq
	\tfrac{\sqrt{M}}{1+\sqrt{M}}\big(1-1/4
	-\sqrt{\tfrac{1}{4\log(M)}} \big) \geq 0.07.
	\end{multline*}
	Note that the constant $\cst[(1)]{\psi, \delta,\gamma}$ does only depend on
	$\psi,\gamma $ and $\delta$, hence 
	it is independent of the parameters $\wSob,\rSob$ and $n$. The claim
	of \cref{theorem:lower_bound} follows by using Markov's inequality,
	which completes the proof.\proEnd
\end{pro}

 \paragraph{Proofs of the lemmata}

\begin{pro}[Proof of \cref{mm:lem:den}]
	For any $h\in C_c^{\infty}(\Rz)$ we can state that $\int_{-\infty}^{\infty} \mathcal S[h](x)dx = [x h(x)]^{\infty}_{-\infty} - \int_{-\infty}^{\infty} h(x) dx=-\int_{-\infty}^{\infty} h(x) dx$ and therefore $\int_{-\infty}^{\infty} \mathcal S^j[h](x)dx = (-1)^j \int_{-\infty}^{\infty} h(x)dx$ for $j\in \Nz$. Thus $\int_{-\infty}^{\infty} \psi_{k,K,\gamma}(x)dx = (-1)^{\gamma}\int_{-\infty}^{\infty} \psi_{k,K}(x)dx =0$ which implies that for any $\delta >0$ and $\bm{\theta}\in\{0,1\}^K$ we have $\int_0^{\infty} f_{\bm{\theta}}(x)dx= 1$.\\
	Now due to the construction \eqref{equation:lobodens} of the functions $\psi_{k,K}$ we easily see that the function  $\psi_{k,K}$ has support on $[1+k/K,1+(k+1)/K]$ which lead to  $\psi_{k,K}$ and $ \psi_{l,K}$ having disjoint supports if $k\neq l$. Here, we want to emphasize that $\mathrm{supp}(\mathcal S[h]) \subseteq \mathrm{supp}(h)$ for all $h\in C_c^{\infty}(\Rz)$. Which implies that $\psi_{k,K,\gamma}$ and $ \psi_{l,K,\gamma}$ have disjoint supports if $k\neq l$, too.
	For $x\in [1,2]^c$ we have $f_{\bm{\theta}}(x)=\exp(-x)\geq 0$. Now let us consider the case $x\in[1,2]$. In fact there is $k_o\in\nsetro{0,K}$ such that $x \in [1+k_o/K,1+ (k_o+1)/K]$ and hence
	\begin{equation*}
	\SoPr[\bm{\theta}](x)= \SoPr[o](x) + \theta_{k_o+1}\delta K^{-\wSob-\gamma} \psi_{k_o,K, \gamma}(x) \geq \exp(-2)  - \delta 2^{\gamma} C_{\gamma, \infty} c_{\gamma}
	\end{equation*}
	since $\|\psi_{k,K,j}\|_{\infty} \leq 2^{j} C_{j, \infty} c_j K^{j}$ for any $k\in \nsetro{0,K}$ and $j\in \Nz$ where $c_j:= \sum_{i=1}^j c_{i,j}$. Now choosing $\delta\leq \delta_o(\psi,\gamma)=\exp(-2)2^{-\gamma} (C_{\gamma,\infty} c_{\gamma})^{-1}$ ensures $f_{\bm{\theta}}(x) \geq 0$ for all $x\in \pRz.$ 
	\proEnd
\end{pro}
	
\begin{pro}[Proof of \cref{lemma:Lag_SobDen}]Our proof starts with the
	observation that for all $t\in \Rz$ we have $\Mel[f_o](t)=\Gamma(1+it)$.
	Now by applying the Stirling formula (see also \cite{BelomestnyGoldenshlugerothers2020} ) we get $|\Gamma(1+it)| \sim |t|^{1/2} \exp(-\pi/2 |t|)$, $|t|\geq 2$, thus for every $s\in \Nz$ there exists $L_{s}$ such that $|f_o|_s^2 \leq L $ for all $L\geq L_s$. \\
	Next we consider $|f_o-f_{\bm{\theta}}|_s$. Let us therefore define first $\Psi_K:= \sum_{k=0}^{K-1} \theta_{k+1} \psi_{k,K}$ and $\Psi_{K,j}:= \mathcal S^j [\Psi_K]$ for an $j\in \Nz$. Then we have
	$|f_o-f_{\bm{\theta}}|_s^2= \delta^2 K^{-2s-2\gamma} |\Psi_{K,\gamma}|_s^2$ where $|\, . \,|_s$ is defined in \eqref{in:eq:MelSob}. Now since for any $j\in \Nz$, it holds that $\mathrm{supp}(\Psi_{K,j}) \subset [1,2], \|\Psi_{K,j}\|_{\infty} < \infty$ we have that $(0,\infty)$ is a subset of the strip of analyticity of $\Psi_{K,j}$. By application of \eqref{equation:Mell_rule} we deduce that $|\Mel[\Psi_{K,s+\gamma}](t)|^2 = (1+t^2)^s|\Mel[\Psi_{K,\gamma}](t)|^2$ and thus
	\begin{align*}
	|\Psi_{K,\gamma}|_s^2&= \int_{-\infty}^{\infty} |\Mel[\Psi_{K,s+\gamma}](t)|^2 dt =2\pi\int_0^{\infty} x |\Psi_{K,s+\gamma}(x)|^2 dx
		\end{align*}
		by the Parseval formula. Since $\psi_{k,K}$ have disjoint support for different values of $k$ we follow that $|\Psi_{k,\gamma}|_s^2 =2\pi \sum_{k=0}^{K-1} \theta_{k+1}^2 \int_0^{\infty} x |\mathcal S^{\gamma+s}[\psi_{k,K}](x)|^2dx$.
		Applying the  Jensen inequality and the fact that $\mathrm{supp}(\psi_{k,K})\subset [1,2]$ leads to
		\begin{align*}
		|\Psi_{k,\gamma}|_s^2&\leq 2 \pi  2^{\gamma+s-1} \sum_{k=0}^{K-1} \sum_{j=1}^{\gamma+s} c_{j,\gamma+s}^2 \int_1^{2} x^{2j+1} K^{2j} \psi^{(j)}(xK-K-k)^2 dx\\
		&\leq 2\pi K^{2(\gamma+s)} 2^{\gamma+s}\sum_{k=0}^{K-1} \sum_{j=1}^{\gamma+s} c_{j,\gamma+s}^2  4^j C_{\psi,s,\gamma}^2 K^{-1} \leq  C_{(\gamma, s)} K^{2(\gamma+s)} 
		\end{align*} 
	Thus $|f_o-f_{\bm{\theta}}|_s^2 \leq C_{( s, \gamma,\delta)}$ and  $|f_{\bm{\theta}}|_s^2 \leq 2(|f_o-f_{\bm{\theta}}|_s^2 + |f_o|_s^2) \leq 2(C_{(s, \gamma,\delta)}+ L_s) =: L_{s, \gamma,\delta}$.
	\proEnd
\end{pro}
\begin{pro}[Proof of \cref{lemma:tsyb_vorb}]  $\ $\\
	Using  that the functions $(\psi_{k,K, \gamma})_{k\in\nsetro{0,K}}$ with different index $k$ have disjoint supports  we get
	\begin{align*}
	\| f_{\bm{\theta}}-f_{\bm{\theta}'}\|_{\omega}^2&= \delta^2 K^{-2s-2\gamma} \| \sum_{k=0}^{K-1} ( \theta_{k+1}-\theta'_{k+1}) \psi_{k,K, \gamma}\|_{\omega}^2 = \delta^2  K^{-2s-2\gamma}\rho(\bm \theta, \bm \theta')  \|\psi_{0,K, \gamma}\|_{\omega}^2
	\end{align*}
	with $\rho(\bm \theta, \bm \theta'):= \sum_{j=0}^{K-1} \1_{\{\bm \theta_{j+1} = \bm \theta'_{j+1}\}}$ the \textit{Hamming distance}. Now the first claim follows by showing that by $ \|\psi_{0,K, \gamma}\|_{\omega}^2 \geq \frac{K^{2\gamma-1}\|\psi^{(\gamma)}\|^2 }{2} $ for $K$ big enough.
	To do so we observe that $ \| \psi_{0,K, \gamma} \|_{\omega}^2 = \sum_{i,j \in \llbracket 1, \gamma\rrbracket} c_{j,\gamma}c_{i,\gamma} \int_0^{\infty} x^{j+i+1} \psi_{0,K}^{(j)}(x)  \psi_{0,K}^{(i)}(x)dx $ and define $\Sigma:= \| \psi_{0,K,\gamma}\|_{\omega}^2 -  \int_0^{\infty}(x^{\gamma} \psi_{0,K}^{(\gamma)}(x))^2 xdx$ 
	\begin{align}\label{equation:low_l2_gamma}
	\| \psi_{0,K, \gamma} \|_{\omega}^2 &= \Sigma + \int_0^{\infty}(x^{\gamma} \psi_{0,K}^{(\gamma)}(x))^2 xdx
	\geq \Sigma + K^{2\gamma-1} \| \psi^{(\gamma)} \|^2 \geq \frac{ K^{2\gamma-1} \| \psi^{(\gamma)} \|^2}{2}
	\end{align}
	as soon as $|\Sigma|\leq \frac{ K^{2\gamma-1} \| \psi^{(\gamma)} \|^2}{2}$. This is obviously true as soon as $K \geq K_{o}(\gamma, \psi)$ and thus $	\| f_{\bm{\theta}}-f_{\bm{\theta}'}\|_{\omega}^2 \geq  \frac{\delta^2\|\psi^{(\gamma)}\|^2 }{2}K^{-2s-1}\rho(\bm{\theta},\bm{\theta}')$ for $K\geq K_o(\psi, \gamma)$.\\
	Now we use the \textit{Varshamov-Gilbert Lemma} (see \cite{Tsybakov2008}) which states that for $K \geq 8$ there existes a subset $\{\bm \theta^{(0)}, \dots, \bm\theta^{(M)}\}$ of $\{0,1\}^K$ with $\bm \theta^{(0)}=(0, \dots, 0)$ such that $\rho(\bm \theta^{(j)}, \bm \theta^{(k)}) \geq K/8$ for all $j ,k\in \llbracket 0, M \rrbracket, j\neq k $ and $M\geq 2^{K/8}$. Applying this leads to $\| f_{\bm\theta^{(j)}}- f_{\bm\theta^{(l)}}\|_{\omega}^2 \geq \frac{\|\psi^{(\gamma)}\|^2\delta^2}{16}  K^{-2s}.
	$\\
	For the second part we have $f_o=f_{\bm\theta^{(0)}}$ and by using $\text{KL}(\widetilde{f}_{\bm \theta}, \widetilde{f}_o) \leq \chi^2(\widetilde{f}_{\bm\theta},\widetilde{f}_o):= \int_{\pRz} (\widetilde{f}_{\bm\theta}(x) -\widetilde{f}_o(x))^2/\widetilde{f}_o(x) dx$ it is sufficient to bound the $\chi$-squared divergence. We notice that $\widetilde{f}_{\bm\theta} -\widetilde{f}_o$ has support in $[0,2]$ since $f_{\bm{\theta}}-f_o$ has support in $[1,2]$ and $g$ has support in $[0,1]$ In fact for $y>2$ holds $\widetilde{f}_{\bm\theta}(y) -\widetilde{f}_o(y)=\int_y^{\infty} (f_{\bm{\theta}}-f_o)(x)x^{-1} g(y/x)dx =0$. Denote further $\Psi_{K,\gamma}:= \sum_{k=0}^{K-1} \theta_{k+1} \psi_{k,K, \gamma}= \mathcal S^{\gamma}[\sum_{k=0}^{K-1} \theta_{k+1} \psi_{k,K}]=: \mathcal S^{\gamma}[\Psi_K]$. Now by using the  compact support property and a single substitution we get
	\begin{align*}
	\chi^2(\widetilde{f}_{\bm\theta},\widetilde{f}_o)
	\leq \widetilde{f}_o(2)^{-1} \|\widetilde{f}_{\bm \theta}- \widetilde{f}_o\|^2 &=\widetilde{f}_o(2)^{-1}  \delta^2 K^{-2s-2\gamma} \| \widetilde{\Psi}_{K,\gamma}\|^2.
	\end{align*}
	 Let us now consider $\| \widetilde{\Psi}_{K,\gamma}\|^2$. In the first step we see by application of the Parseval that $\|\widetilde{\Psi}_{K,\gamma}\|^2 = \frac{1}{2\pi} \int_{-\infty}^{\infty} |\Mela{\widetilde{\Psi}_{K,y}}{1/2}(t)|^2dt$. Now for  $t\in \Rz$, we see by using the multiplication theorem for Mellin transforms  that $\Mela{\widetilde{\Psi}_{K,\gamma}}{1/2}(t)= \Mela g {1/2}(t) \cdot \Mela{\mathcal S^{\gamma}[\Psi_{K}]}{1/2}(t)$.  Again we have $\Mela{ \mathcal S^{\gamma}[\Psi_{K}]}{1/2}(t)=(1/2+it)^{\gamma} \Mela{ \Psi_{K}}{1/2}(t)$. Together with assumption \textbf{[G1']} we get
	\begin{align*}
	\|\widetilde{\Psi}_{K,\gamma}\|^2  \leq \frac{C_1(g)}{2\pi}\int_{-\infty}^{\infty} |\Mela{ \Psi_{K}}{1/2}(t)|^2 dt= C_1(g) \| \Psi_K\|^2  \leq C_1(g) \|\psi\|^2. 
	\end{align*}
Since $M \geq 2^K$ we have thus
	$\text{KL}(\widetilde{f}_{\bm\theta^{(j)}},\widetilde{f}_{\bm\theta^{(0)}})\leq \frac{C_1(g)\|\psi\|^2}{\widetilde f_o(2)\log(2)} \delta^2 \log(M) K^{-2s-2\gamma-1}.$ \proEnd
\end{pro}

%%% Local Variables:
%%% mode: latex
%%% TeX-master: "_0SCORDEMME"
%%% End:

% --------------------------------------------------------------------
% <<BibFile>>
% --------------------------------------------------------------------
\bibliography{SCORDEMME} 
\end{document}